\documentclass[a4paper,12pt]{article}
\usepackage[utf8]{inputenc}
\usepackage[T1]{fontenc}
\usepackage{palatino, mathpazo,amsthm}
\usepackage{url}
\usepackage{hyperref}
\usepackage[margin=2.5cm,nohead]{geometry}
\usepackage{parskip}
\usepackage{xcolor}
\usepackage{graphicx,color}
\usepackage{boxedminipage}

\usepackage{epsfig}

{\begin{boxedminipage}{#1}}
{\end{boxedminipage}}%

\usepackage{rasmus2,natbib}

\hypersetup{
  colorlinks,
  citecolor=black,%
  filecolor=black,%
  linkcolor=black,%
  urlcolor=black
}

\graphicspath{{.}{./code/pix-auto/}}

\title{Quick inference for log Gaussian Cox processes with non-stationary underlying random fields}

\author{Ji\v r\' i Dvo\v r\' ak\footnote{Department of Probability and Mathematical Statistics, Faculty of Mathematics and Physics, Charles University, Prague, Czech Republic; email: dvorak@karlin.mff.cuni.cz}, 
Jesper M\o ller\footnote{Department of Mathematical Sciences, Aalborg University, Denmark; email: jm@math.aau.dk},
Tom\' a\v s Mrkvi\v cka\footnote{Department of Applied Mathematics and Informatics, Faculty of Economics, University of South Bohemia, \v{C}esk\'e Bud\v{e}jovice, Czech Republic; email: mrkvicka.toma@gmail.com},  
\ and Samuel Soubeyrand \footnote{BioSP, INRA, 84914, Avignon, France; email: samuel.soubeyrand@inra.fr}}
\date{July 31, 2019}

\begin{document}
\maketitle

{\bf Abstract:} 
For point patterns observed {\it in natura}, spatial heterogeneity is more the rule than the exception. In numerous applications, this can be mathematically handled by the flexible class of log Gaussian Cox processes (LGCPs); in brief, a LGCP is a Cox process driven by an underlying log Gaussian random field (log GRF). This allows the representation of point aggregation, point vacuum and intermediate situations, with more or less rapid transitions between these different states depending on the properties of GRF. Very often, the covariance function of the GRF is assumed to be stationary. In this article, we give two examples where the sizes (that is, the number of points) and the spatial extents of point clusters are allowed to vary in space. To tackle such features, we propose parametric and semiparametric models of non-stationary LGCPs where the non-stationarity is included in both the mean function and the covariance function of the GRF. 
Thus, in contrast to most other work on inhomogeneous LGCPs, second-order intensity-reweighted stationarity is not satisfied and the usual two step procedure for parameter estimation based on e.g.\ composite likelihood does not easily apply. Instead 
we propose a fast three step procedure based on composite likelihood. We apply our modelling and estimation framework to analyse datasets dealing with fish aggregation in a reservoir and with dispersal of biological particles.



\section{Introduction}

\subsection{Motivation and objective}\label{s:motivation}

Inhomogeneity and aggregation (also called clumping or clustering) in spatial point patterns are frequently observed in nature because of spatially-structured covariates and random effects. Being able to infer these characteristics generally contributes to unravel influential underlying processes.

\begin{figure}[t]
\begin{center}
\includegraphics[height=\textwidth,angle=-90]{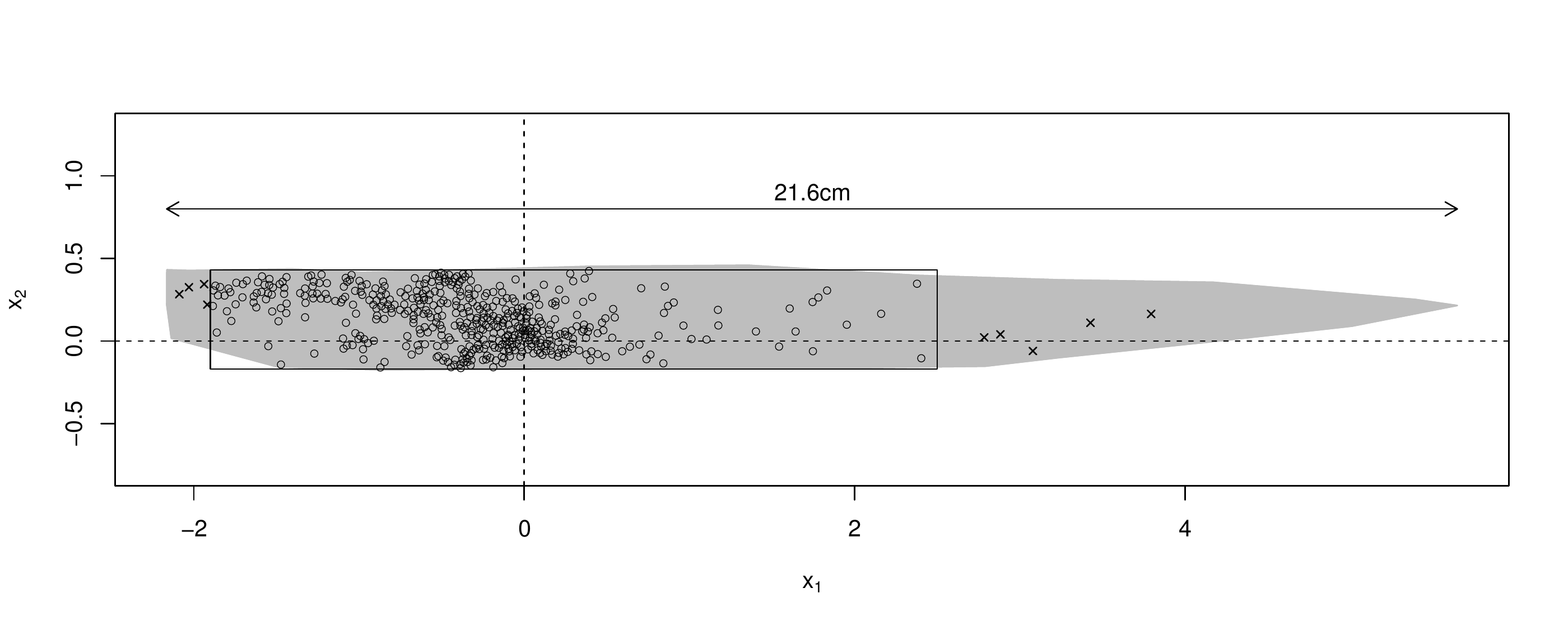} 
\end{center}
\caption{Locations (open circles and crosses) of fungal lesions on a plant leaf (grey shape). The fungal lesions were formed after the dispersal of spores emitted from one lesion located at the intersection of the two dashed lines. The estimation algorithm was applied to the point pattern within the rectangle (open circles; 476 points), whereas locations of lesions outside the rectangle (crosses; 9 points) were not used.}
 \label{f:particles}
\end{figure}

\begin{figure}[t]
\centering 
  \includegraphics[width=\textwidth]{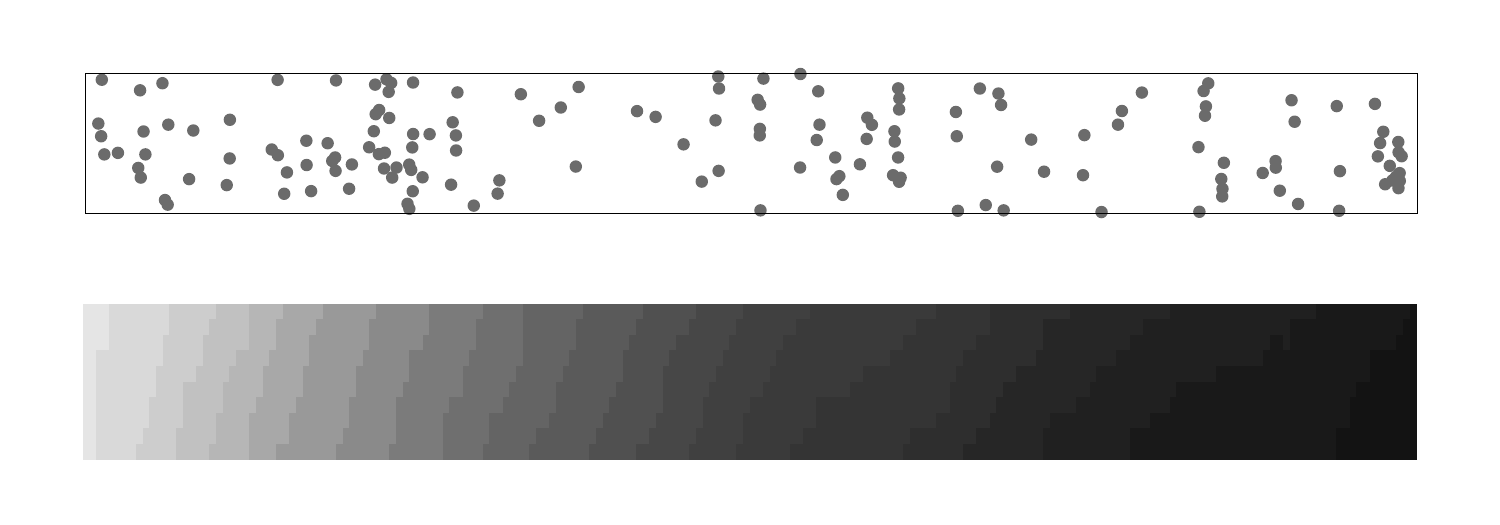} 
   \caption{Top panel: positions of fish recorded along the boat trace; only a 200 meter long part of the data is shown, the width of the observation window is 10 meters. Bottom panel: water depth; values range from 5 (dark grey) to 29 meters (light grey).}\label{f:fishdata}
\end{figure}

For example, Figure~\ref{f:particles} shows the locations of 485 brown rust lesions on a wheat leaf. The lesions were obtained by the dispersal of spores emitted from one mother lesion which we refer to as the point source (the point (0,0) in the figure; further details are given in Section~\ref{sec:particles}). Our aim is to describe how the intensity and the random clumping of lesions depends on the distance to the point source. 
Another example is given in Figure~\ref{f:fishdata}, which shows a part of a dataset consisting of the positions of 706 fish recorded along a boat trace. Here, our aim is to describe how the covariate shown in the figure (the spatially varying depth of water) influences the random organisation of the fish shoals (further details are given in Section~\ref{sec:fish}).

Inhomogeneity and aggregation in spatial point patterns are usually described in terms of the intensity function $\rho$ and the pair correlation function $g$ defined as follows.  
We consider an observed planar point pattern (such as in Figures~\ref{f:particles}--\ref{f:fishdata}) as a realization of the restriction of a planar point process $X$ to a bounded observation window $W\subset\mathbb R^2$. 
Excluding the case of multiple points, we can view $X$ as a locally finite random subset of $\mathbb R^2$. The points in $X$ are  called {\it events} the points in $X$, and we denote by $N(A)=\#(X\cap A)$ the number of events within a region $A\subseteq \mathbb R^2$ (strictly speaking $A$ needs to be a Borel set, but we suppress measure theoretical details; Section~\ref{s:background} provides more details on planar point processes).
Now, it is 
assumed that whenever $A\subset\mathbb R^2$ is bounded, the expected number of events in $A$ is finite and given by
\begin{equation}\label{e:rhodef}
\mathrm E\{N(A)\}=\int_A\rho(u)\,\mathrm du,
\end{equation}
and whenever $A,B\subset\mathbb R^2$ are bounded and disjoint,
the covariance between $N(A)$ and $N(B)$ exists and is given by
\begin{equation}\label{e:gdef}
\mathrm{Cov}\{N(A),N(B)\}=\int_A\int_B\rho(u)\rho(v)\{g(u,v)-1\}\,\mathrm du\mathrm dv.
\end{equation}
Whilst the intensity function describes the inhomogeneity, the pair correlation function describes the positive ($g>1$) or negative ($g<1$) dependence between pairs of counts; generally speaking, $g>1$ corresponds to aggregation of the events.

Thus, with a real dataset such as those shown in Figures~\ref{f:particles}-\ref{f:fishdata}, inhomogeneity and aggregation can be described by estimating $\rho$ and $g$. This estimation can be tackled by constructing flexible parametric models for $\rho$ and $g$ and evaluating the most {\it likely} values for the model parameters. In this article, $\rho$ and $g$ are modelled in the framework of parametric and semiparametric models for spatial Cox processes, where we include first- and second-order non-stationarity into $\rho$ and $g$, respectively. We develop a fast method for parameter estimation based on only $\rho$ and $g$, which enables processing  large point pattern data sets.

\subsection{Spatial Cox processes}\label{s:Coxproc}   

Recall that $X$ is a planar Poisson process with intensity function $\rho$ if 
for any bounded region $A\subset\mathbb R^2$, $N(A)$ is Poisson distributed with mean $\int_A\rho(u)\,\mathrm du$, and conditioned on $N(A)=n$, the $n$ events in $A$ are independent and identically distributed with density proportional to $\rho$ restricted to $A$.
Poisson processes provide models for first order inhomogeneity as specified by $\rho$, whilst they are not flexible enough for modelling second and higher order spatial dependence properties as  reflected, for instance, by the pair correlation function; see e.g. the discussion in \cite{moeller:waagepetersen:07} and the references therein. 

In contrast, Cox processes \citep{cox:55} provide a large and flexible class of models for aggregation, where expressions for $\rho$ and $g$ are available. Recall that $X$ is a planar Cox process driven by a random field $\Lambda:\mathbb R^2\mapsto[0,\infty)$ if $X$ conditioned on $\Lambda$ is a Poisson process with intensity function $\Lambda$. For distinct points $u,v\in\mathbb R^2$, we have 
\[\rho(u)=\mathrm E\{\Lambda(u)\},\qquad\rho(u)\rho(v)g(u,v)=\mathrm E\{\Lambda(u)\Lambda(v)\},\]
 provided these functions are locally integrable, that is, the integrals in \eqref{e:rhodef}--\eqref{e:gdef} are well-defined. Usually, for specific parametric Cox process models (including the models in this paper), we have positive dependence ($g>1$), 
although it is possible to construct Cox processes where $g(u,v)<1$ for some locations $u,v\in\mathbb R^2$. 
\cite{moeller:waagepetersen:07} suggested separating first order inhomogeneity specified by $\rho$ from random effects specified 
 by a residual random field $\Lambda_0$ such that for all $u\in\mathbb R^2$,
\begin{equation}\label{e:inhomCox}
\Lambda(u)=\rho(u)\Lambda_0(u),\qquad \mathrm E\{\Lambda_0(u)\}=1.
\end{equation} 
This allows us
to model independently the intensity and the kind of aggregation 
which is not explained by the intensity.
Typically, when covariate functions $z_1,\ldots,z_I:W\mapsto\mathbb R$ are available,
a log linear intensity function $\xi:=\log\rho$ is assumed:
\begin{equation}\label{e:log-linear-rho}
\xi(u)=\sum_{i=0}^I z_i(u)\beta_i
\end{equation}
where $z_0\equiv 1$ corresponds to an intercept and $\beta_0,\ldots,\beta_I$ are real regression parameters.
Note that the Cox process driven by $\Lambda_0$ has also pair correlation function $g$.
In \cite{moeller:waagepetersen:07} and subsequent work (e.g.\ \cite{diggle:13} and the references therein), $\Lambda_0$ is assumed to be stationary, that is, its distribution is invariant under translations in $\mathbb R^2$; this implies that $g(u,v)=g_0(u-v)$ depends only on $u-v$. 
However, we shall not make this assumption as we will need more flexibility when modelling $\Lambda_0$ for datasets as in 
 Figures~\ref{f:particles}--\ref{f:fishdata}.

\subsection{Spatial Cox processes with non-stationary underlying random fields}\label{sec:CCP}

In this paper, we do not assume that $\Lambda_0$ is stationary, because we want to allow cluster sizes and extents to depend on spatial location. This point of view was also taken in \cite{tomas:14} and
\cite{tomas:samuel:16} who considered the detection of different types of inhomogeneities by using 
 shot noise Cox process 
 models, as  discussed further in Section~\ref{s:LGCP}. Here, in contrast, we will adopt the Cox process framework, but will extend it to particularly flexible cases.

Consider a log Gaussian Cox process (LGCP), that is, a Cox process as above when the log residual random field $Y_0=\log\Lambda_0$ is a Gaussian random field (GRF) such that $\mathrm E[\exp\{Y_0(u)\}]=1$. Thus, denoting the 
covariance function of the GRF by $c(u,v)=\mathrm{Cov}\{Y_0(u),Y_0(v)\}$ and the variance function by $\sigma^2(u)=c(u,u)$, the GRF has mean function $\mathrm E\{Y_0(u)\}=-\sigma^2(u)/2$. Often (including the cases of this paper) $c\ge0$ or equivalently $g\ge1$ (as discussed further in Section~\ref{s:LGCP}). We will consider models where $\rho(u)$ describes the overall `average' inhomogeneity, and  
$\sigma(u)=\sqrt{\sigma(u)^2}$ controls the spatial varying weights (in the sense of the number of events) of the clumps/clusters. Additionally, inspired by the space transformation approach in \cite{sampson:guttorp:92} and subsequent work \citep[e.g.][]{perrin2000,fouedjio2015}, we also allow the covariance function to be of the form
\begin{equation}\label{e:cov-model}
c(u,v)=\sigma(u)r(\|\varphi(u)-\varphi(v)\|)\sigma(v).
\end{equation}
Here, $r:[0,\infty)\mapsto[-1,1]$ is a positive semi-definite function with $r(0)=1$, $\|\cdot\|$ denotes the Euclidean distance in $\mathbb R^2$, and $\varphi:\mathbb R^2\mapsto\mathbb R^2$ is a bijective     
mapping. 
The specification of $\varphi$ and $r$ allows us to model the spatial varying shape and scale of clusters. For instance, we can model elliptical shapes by linear functions $\varphi(u)=Au$ \citep{sampson:guttorp:92,moeller:toftager:14}, and model a varying spatial scale with $\varphi(u)=\|u\|^{\delta-1}Au$ or 
$\varphi(u)=\log(1+\delta\|u\|)\|u\|^{-1}Au$, 
where $A$ is a $2\times2$ regular matrix and $\delta>0$ ; the latter is illustrated later in Section~\ref{sec:particles} when $A$ is the identity matrix. Moreover, the geostatistical literature \citep{cressie1993,chiles1999,stein1999} suggests numerous parametric 
models for the function $r$ which corresponds to the correlation function for a stationary and isotropic random field.

In summary, we let $Q$ be a  zero-mean and unit-variance GRF with an isotropic correlation function $r(t)=\mathrm{Cov}\{Q(u),Q(v)\}$ for $t=\|u-v\|$, and consider 
\begin{equation}\label{e:logL}
\log\Lambda(u) 
=\xi(u)+Y_0(u),\qquad Y_0(u)=-\sigma^2(u)/2+\sigma(u)Q\{\varphi(u)\},\qquad u\in\mathbb R^2,
\end{equation}
where $\xi$ is given by \eqref{e:log-linear-rho}. Examples of realizations of LGCPs based on such random fields are shown in Figures~\ref{fig:fish.simul} and \ref{fig:particles.simul}. To the best of our knowledge, such LGCPs have only been studied in the literature when $\sigma$ is constant and $\varphi$ is the identity mapping.

\subsection{Outline and software}

The remainder of this paper is organized as follows. Section~\ref{s:background} provides the needed background material on Cox processes and parametric inference procedures, and summarizes the attractive properties of LGCPs. Section~\ref{sec:3step-estim-check} presents our proposal for performing rapid estimation of parameters of a LGCPs with non-stationary underlying random field using a three step procedure based on composite likelihoods.
Our approach is illustrated in Section~\ref{sec:fish} in a case study concerning fish aggregation (Figure~\ref{f:fishdata}), and in Section~\ref{sec:particles} in a  case study dealing with the dispersal of particles from a point source (Figure~\ref{f:particles}). Section~\ref{s:conclusion} summarizes our findings in Sections~\ref{sec:3step-estim-check}--\ref{sec:particles}. Finally, Appendices A--C provide details related to the data in Figure~\ref{f:particles} and the space transformation function $\varphi$.

Computations were carried out with the {\bf R} statistical software, including the packages {\tt spatstat}, {\tt RandomFields} and {\tt GET}. Codes for running the estimation procedure are available upon request.

\section{Background}\label{s:background}

Henceforth, we assume that $X$ is a Cox process driven by a random intensity function $\Lambda$ as in \eqref{e:inhomCox} and a parametric statistical model has been specified in terms of an unknown parameter vector $\theta=(\beta,\nu)$, where $\beta$ is used to parametrize $\rho(u)=\rho(u;\beta)$, $\nu$ to parametrize the distribution of $\Lambda_0$, and $\beta$ and $\nu$ are variation independent. Thus, $g(u,v)=g(u,v;\nu)$ depends only on $\nu$.

\subsection{Composite likelihood estimation}\label{subsec:CLE}
 
Suppose $W\subset\mathbb R^2$ is a bounded observation window and a realization $X\cap W=\{x_1,\ldots,x_n\}$ has been observed. Due to the computational obstacles (especially if $n$ is large) 
related to the calculation of the likelihood function, computationally easier approaches for statistical inference based on functional summary statistics, including $\rho(u;\beta)$ and $g(u,v;\nu)$, together with estimation equations (obtained e.g.\ by minimum contrast methods, composite likelihoods or Palm likelihoods) have been suggested, see the review in \cite{moeller:waagepetersen:16}. 
We will concentrate on the first and second order composite likelihoods \citep{guan:06,waagepetersen:07,tanaka:ogata:stoyan:08,waagepetersen:guan:09,guan:jalilian:waagepetersen:15,zhuang:15}.
The first-order log composite likelihood is given by
\begin{equation}\label{e:CL1}
CL_1(\beta)=\sum_{i=1}^n\log\rho(x_i;\beta)-\int_W\rho(u;\beta)\,\mathrm du.
\end{equation}
This is also known as the Poisson log likelihood 
and it can be interpreted as the composite log likelihood 
 for binary indicators $I_k=1\left(N(C_k)>0\right)$  of the presence of events in cells $C_k$ of an infinitesimal partition of $W$ \citep{schoenberg:05,waagepetersen:07,guan:loh:07,moeller:waagepetersen:07}. Further, for a tuning parameter $R>0$, the second-order log composite likelihood is given by
\begin{align}
CL_2(\beta,\nu) 
=& \sum_{1\le i<j\le n}1\left(\|x_i-x_j\|\le R\right)\Big[\log\left\{\rho\left(x_i;\beta\right)\rho\left(x_j;\beta\right)g\left(x_i,x_j;\nu\right)\right\}\nonumber  \\
&
 -\log \int_W\int_W 1\left(\|u-v\|\le R\right)\rho\left(u;\beta\right)\rho\left(v;\beta\right)g\left(u,v;\nu\right)\,\mathrm du\,\mathrm dv  \Big]
\label{e:CL2}
\end{align}
and it can be interpreted as the log composite likelihood for the binary indicators $I_kI_l$ of simultaneous occurrence of events in $R$-close pairs of distinct cells $C_k$ and $C_l$ \citep{waagepetersen:07,moeller:waagepetersen:07}.

The following two step procedure is widely used for finding estimates $(\hat\beta,\hat\nu)$ of $(\beta,\nu)$ \citep{waagepetersen:guan:09,prokesova:dvorak:jensen:15, M12}. First, $\hat\beta$ is obtained by maximizing the first-order log composite likelihood \eqref{e:CL1}; the log-linear form  for $\rho$ (see Equation \eqref{e:log-linear-rho}) allows a fast estimation of $\hat\beta$ by using {\tt spatstat}. Second, $\hat\nu$ is obtained by maximizing the second-order log composite likelihood $CL_2(\hat\beta,\nu)$.

\subsection{Pros of using log Gaussian Cox processes}\label{s:LGCP}

The two main classes of Cox processes are log Gaussian Cox processes (LGCPs; see Section~\ref{sec:CCP}) and  shot noise Cox processes (SNCPs; see \cite{moeller:03}, \cite{moeller:waagepetersen:03}, and the references therein). For a SNCP $X$ driven by \eqref{e:inhomCox}, the residual random field is of the form:
\[\Lambda_0(u)=\sum_{c\in C}k_c(u),\]
where $C$ is a planar Poisson process with intensity function $\rho_C$, and where $k_c$ is a non-negative function such that $\mathrm E\{\Lambda_0(u)\}=1$ or equivalently
\begin{equation}\label{e:int1}
\int\rho_C(c)k_c(u)\,\mathrm dc=1,\qquad u\in\mathbb R^2.
\end{equation} 
Assuming $\rho_c(u)=\rho(u)k_c(u)$ is integrable, then $X$ conditioned on $C$ is distributed as $\cup_{c\in C}X_c$, where
$X_c$ is a finite Poisson process with intensity function $\rho_c$, and where these Poisson processes for all $c\in C$ are independent.  Therefore $C$ is called the centre process and $X_c$ the cluster with centre $c\in C$. Moreover,
\begin{equation}\label{e:int2}
g(u,v)=1+\int k_c(u)k_c(v)\rho_C(c)\,\mathrm dc,\qquad u,v\in\mathbb R^2,\ u\not=v,
\end{equation}
provided the integral exits, and then $g\ge 1$. In general, unless $\Lambda_0$ is stationary or equivalently $\rho_C(u)$ is constant 
and $k_c(u)=k_0(u-c)$ for all $c,u\in\mathbb R^2$, the integrals in \eqref{e:int1}--\eqref{e:int2} are not expressible on closed form; and even more complicated integrals appear for the third and higher order joint intensities. For example, \cite{tomas:samuel:16}
considered non-stationary SNCPs, for which they needed to evaluate the integrals by numerical methods, and they proposed a Bayesian approach based on a MCMC algorithm.
On one hand, their approach provides a detailed analysis of uncertainty and, on the other hand, it requires very time-consuming calculations.   

In the aim of proposing fast composite likelihood based estimation procedure for second-order non-stationary Cox processes, we find LGCPs as given by \eqref{e:inhomCox} and \eqref{e:logL} attractive for the following reasons. 
First, the pair correlation function is in a simple one-to-one correspondence to the covariance function $c$ of the Gaussian process $\log\Lambda_0$: 
\[g(u,v)=\exp\left\{c(u,v)\right\}.\]
Thus, the distribution of $X$ is uniquely determined by $\rho$ and $g$. This is in line with the first and second order composite likelihood approach. 
 Second, 
edge effects are not playing a role as the distribution of $X$ restricted to a bounded observation window $W$ is determined by $\rho$  restricted to $W$ together with
the distribution of the Gaussian process $\log\Lambda_0$ restricted to $W$.
Finally, a LGCP possesses the following attractive properties, although these will not be exploited in this article: the reduced Palm distributions 
of a LGCP are also LGCPs with unchanged pair correlation function $g$ \citep{coeurjolly:moeller:waagepetersen:LGCP}, and third and higher order intensities are easily expressed in terms of $\rho$ and $g$ \citep{moeller:syversveen:waagepetersen:98}.

\section{Parameter estimation and model check for Cox processes driven by non-stationary random fields}\label{sec:3step-estim-check}

Throughout this section, for specificity and simplicity, we assume that $X$ is a LGCP governed by \eqref{e:inhomCox} and \eqref{e:logL}. However, the description of the estimation procedure in Section~\ref{sec:3-steps} can easily be modified to other cases of parametric models for Cox processes driven by non-stationary random fields, including the case of SNCPs. 

\subsection{A three step composite likelihood based estimation approach}\label{sec:3-steps}

For parametric models of LGCPs 
with the non-stationary covariance structure introduced in Section \ref{sec:CCP}, particularly in connection to the analysis of the data presented in Sections~\ref{s:Tomasdata} and \ref{s:Samueldata},
we noticed after various experiments, some of which are discussed in~\ref{app:comparison}, that the two step procedure presented above yields strongly biased estimates, in particular for the parameters governing the second-order non-stationarity. This  relates to the fact that a nonlinear function $\varphi$ makes it impossible to choose a single value of $R$ that would be appropriate in the whole observation window. In other words, the model is too complex (compared to the size of datasets we consider) to enable reliable estimation of interaction parameters in a single step as performed in the two step estimation procedure.
In this section, we therefore suggest a three step procedure based on composite likelihoods that circumvents the non-stationarity issue as detailed below. 

We split the dataset with respect to a spatial covariate function, say $f:\mathbb R^2\to\mathbb R$, which will be used to delineate disjoint areas of $W$ over which the point process $X$ is approximately second-order stationary. The three step procedure has to be adapted to the model of interest and, here, we only provide its skeleton.  Let $f_0<\ldots<f_K$ be a sequence of real numbers,  
and let $W_k=\{u\in W:f_{k-1}\le f(u) < f_k\}$ for $k=1,\ldots,K$. We have in mind that $f_0,\ldots,f_K$ are chosen such that (i) $X$ is approximately second-order stationary over $W_k$, (ii) $X\cap W_1,\ldots,X\cap W_K$ contain approximately the same numbers of points and (iii) $W=W_1\cup\cdots\cup W_K$. The division of the dataset to smaller parts which are assumed to be second-order stationary, which is performed in the proposed procedure, makes the estimation more stable as it can be seen in~\ref{app:comparison}. Due to the assumption of second-order stationarity in individual subwindows we can expect that the proposed procedure will be successful in cases of smoothly changing second-order structure.

Now, define the approximate second-order log composite likelihood restricted to $W_k$ and based on a second-order stationary assumption by
\begin{align}
CL_2^k(\beta,\nu_0) 
=& \sum_{\begin{smallmatrix}1\le i<j\le n\\ x_i,x_k\in X\cap W_k\end{smallmatrix}}1\left(\|x_i-x_j\|\le R_k\right)\Big[\log\left\{\rho\left(x_i;\beta\right)\rho\left(x_j;\beta\right)g_0\left(x_i,x_j;\nu_0\right)\right\}\nonumber  \\
&
 -\log \int_{W_k}\int_{W_k} 1\left(\|u-v\|\le R_k\right)\rho\left(u;\beta\right)\rho\left(v;\beta\right)g_0\left(u,v;\nu_0\right)\,\mathrm du\,\mathrm dv  \Big].
\label{e:CL2.split}
\end{align}
Here, each $R_k>0$ is a tuning parameter that may depend on $k$, $g_0$ is the pair correlation function obtained by assuming that the covariance structure is stationary (that is, $\sigma$ is constant and $\varphi(u)=u$), and $\nu_0$ is the corresponding set of parameters (typically the constant value of $\sigma$ and the range parameter arising in the correlation function $r$).
In addition, suppose the parameter $\nu$ related to the covariance structure is a vector
$\nu=(\nu^\text{split},\nu^\text{joint})$, where $\nu^\text{split}$ contains the components of $\nu$ that will be estimated via the approximate second-order composite likelihood applied to the restricted point processes $X\cap W_1,\ldots,X\cap W_K$, and $\nu^\text{joint}$ contains the components of $\nu$ that will be estimated via the second-order composite likelihood applied to the whole point process $X$.

Then the estimation procedure consists of the three following steps.
\begin{enumerate}
\item If the first-order intensity $\rho$ is modelled by a function parameterized by $\beta$ (like the log-linear function in Equation \eqref{e:log-linear-rho}), $\hat\beta$ is obtained by maximizing the first-order log composite likelihood \eqref{e:CL1}:
$$\hat\beta=\text{argmax}_{\beta}CL_1(\beta).$$
Otherwise, a kernel-smoothing approach is used to estimate $\rho$.
\item The part $\nu^\text{split}$ of $\nu$ is estimated by firstly estimating $\nu_0$ for each $k$ via the maximization of the approximate second-order log composite likelihood \eqref{e:CL2.split}:
$$\hat\nu_{0,k}=\text{argmax}_{\nu_0}CL_2^k(\hat\beta,\nu_0);$$
and by secondly fitting a regression with $\hat\nu_{0,k}$ as the response variable and the average value of $f$ in $W_k$, namely $\bar F_k=\frac{1}{\#{X\cap W_k}}\sum_{x\in X\cap W_k} f(x_k)$, as the explanatory variable to estimate each component of $\nu^\text{split}$. In Sections~\ref{sec:fish} and \ref{sec:particles}, we show how these regressions are implemented in practice.
\item After the transformation of the point pattern and the observation window by $\varphi$, the remaining part $\nu^\text{joint}$ of $\nu$ is estimated by maximizing the second-order log composite likelihood \eqref{e:CL2} with respect to $\nu$ and by ignoring the estimate obtained for $\nu^\text{split}$:
$$\hat\nu^\text{joint}=\text{argmax}_{\nu}CL_2(\hat\beta,\nu).$$
\end{enumerate}

Remark 1: If no space deformation is assumed in the model, steps 2 and 3 can be reversed and, in this case, the estimate of $\nu^\text{joint}$ might be plugged in for estimating $\nu^\text{split}$. 

Remark 2: When applying the presented methodology to a new dataset the user has to make several choices in the modelling step. The case studies in Sections~\ref{sec:fish} and \ref{sec:particles} can provide some guidance. Generally speaking, either a parametric model is assumed for the intensity function (see Section~\ref{sec:particles}) or it is estimated nonparametrically (Section~\ref{sec:fish}); either the autocorrelation function of the GRF is considered to be translation invariant (Section~\ref{sec:fish}) or not (Section~\ref{sec:particles}); either the variance of the GRF is treated as being constant (Section~\ref{sec:particles}) or not (Section~\ref{sec:fish}).

Remark 3: A simulation study assessing the estimation efficiency is generally required to select the components of $\nu$ appearing in $\nu^\text{joint}$ and $\nu^\text{split}$, and to select the tuning parameters $K$, $R_k$ and $f_k$.

Remark 4: Using the second-order composite likelihood as if the point pattern is second-order stationary and isotropic provides a computational advantage: the double integral in the composite likelihood criterion \eqref{e:CL2.split} is then of the form
$$\int_{W_k}\int_{W_k} 1\left(\|u-v\|\le R_k\right)h(u,v)g_0\left(\|u-v\|;\nu_0\right)\,\mathrm du\,\mathrm dv$$
and can be computed as the Lebesgue-Stieltjes integral
\begin{align}
  \int_0^{R_k} g_0(t;\nu_0)\,\mathrm dH_k(t) \label{eq:LSint}
\end{align}
where  $H_k(t)=\int_{W_k}\int_{W_k} 1\left(\|u-v\|\le t\right)h(u,v)\,\mathrm du\,\mathrm dv$. Note that $H_k(t)$ does not depend on the parameter $\nu_0$. Thus,  $H_k(t)$ is calculated only once and only the integral $\int_0^{R_k} g_0(t;\nu_0)\,\mathrm dH_k(t)$ is computed in each step of the iterative optimization of the composite likelihood criterion.

\subsection{Test based on global envelopes of the inhomogeneous $J$-function}\label{sec:model-check}

Our fitted models were checked with plots of global envelopes and accompanying tests \citep{MyllymakiEtal2017,MrkvickaEtal2017}. The global envelopes can be  based on various functional summaries of point processes, e.g.\ topological characteristics, morphological characteristics, and third-order characteristics. These statistics can be viewed as appropriate for model checking because they are not used in the estimation procedure considered above. However, they do not directly reflect the dependence of the point pattern structure on spatial covariates, although such a dependence should be taken into account in model checking of point process models with a spatially changing structure.  

To reveal the dependence on a given spatial covariate in connection to a LGCP, we have used the inhomogeneous $J$-function \citep{lieshout:11} computed on several subwindows, which are determined according to the spatial covariate. For example, in the group dispersal model, the subwindows are generated with respect to the distance from the source point emitting the particles. The first subwindow is the inner disk centred around the source point, the second subwindow is the subsequent circular ring, and so on. Note that only the point lying in each subwindow were used for the computation of each inhomogeneous $J$-function.

\section{Fish aggregation}\label{sec:fish}

\subsection{Data and model}\label{s:Tomasdata}

The point pattern studied in this section is part of a larger dataset studied in \cite{M12} and \cite{muska:etal:16}. It consists of the locations of 706 fish observed during a summer day
 along the trace made by a boat in the middle part of the inland water reservoir Rimov in Czech Republic. 
The trace is 916 meter long and the boat crossed the river twice.
The fish locations were recorded within a distance of 10 to 20 meters from the boat and in a
depth of 1 to 1.75 meters. The data we analysed are given by the $(x,y)$-coordinates of the 3D-locations of fish, where the
observation window $W$ is of size $916\times 10$ meters (note that $W$ is larger than the region shown in Figure~\ref{f:fishdata}).
In fact, as the reservoir is thermally stratified, the majority of fish
are in the upper 5 meters of the water column during the summer, and the horizontal distance is considered to be more important than the vertical distance between fish \citep{Jarolim}. 
The aim is to study the effect of a spatial covariate, namely the water depth denoted by $z(u)$, $u\in W$, on the distribution of fish shoals and, more exactly, on the dispersion of this distribution. 

We used a semiparametric LGCP, with an intensity function estimated in a non-para\-metric way (as detailed below) and where the ingredients of the covariance function in \eqref{e:cov-model} is specified as follows.
The standard-deviation depends linearly on the spatial covariate:
\begin{equation}\label{e:sigma-tomas}
 \sigma(u) = \gamma_0 + \gamma_1 z(u),\qquad u\in W,
\end{equation}
where $\gamma_0$ and $\gamma_1$ are real parameters such that $\sigma(u)>0$ for all $u\in W$. 
Further, there is no space transformation, that is, $\varphi$ is the identity function. Finally, the  
 auto-correlation function $r$ is exponential with parameter $\alpha>0$: 
\begin{equation*}\label{e:expcorr2}
r(t)=\exp\left(- \alpha t\right),\qquad t\ge0.
\end{equation*}
This LGCP can generate point patterns  
corresponding to non-stationary aggregation of fish, with plenty of small shoals or a few large shoals or something in between depending on the covariate values and the parameter values.

\subsection{Implementation of the three step procedure}\label{s:parest-tomas}

We adapted the three step estimation procedure of Section \ref{sec:3-steps}, with
$\nu^\text{split}=(\gamma_0,\gamma_1)$ and $\nu^\text{joint}=\alpha$, and with steps as follows.

Step 1: For the estimation of the intensity function $\rho$, 
we used a kernel smoothing method with a large bandwidth (Gaussian kernel with standard deviation equal to 50 meters) in order to allow the estimated $\rho$ function to vary along the boat trace.

Step 2: For the estimation of $\gamma_0$ and $\gamma_1$, 
we started by letting $f\equiv z$, and defined $K=5$ disjoint subwindows $W_1, \ldots, W_5$ as areas corresponding to different ranges of water depth as determined by the values of $f_0,\ldots,f_5$. We chose these values such that the subwindows contain approximately the same number of points, e.g.\ $W_1$ contains approximately 141 points lying in the most shallow water. The value of $K$ was chosen as a compromise between the requirement that $f$ is approximately constant on each $W_k$ and the necessity to have a reasonable number of points in each $W_k$ to allow for reliable estimation. The choices are rather subjective as is the case also for other estimation procedures based on estimating equations and therefore some experimentation is recommended. Then, $\nu_{0,k}=(\sigma_{0,k},\alpha_{0,k})$ was estimated by maximizing the approximate second-order log composite likelihoods based on $X\cap W_k$. Let $\hat\sigma_{0,k}$ and $\hat\alpha_{0,k}$ denote the estimates of $\sigma_{0,k}$ and $\alpha_{0,k}$, respectively, and let $\bar F_k$ denote the average value of $f$ over $W_k$. By assuming that the function $f$ is approximately constant (and approximately equal to $\bar F_k$) over $W_k$, the standard-deviation function $\sigma$ can be approximated by
\begin{align*}
\sigma(u)=\gamma_0+\gamma_1 z(u)\approx \gamma_0+\gamma_1 \bar F_k\qquad\mbox{for }u\in W_k.
\end{align*}
Consequently, we used $\hat\sigma_{0,k}$ as a natural estimate of $\gamma_0+\gamma_1 \bar F_k$. 
Finally, the linear regression model $\hat\sigma_{0,k}=\gamma_0+\gamma_1 \bar F_k+\epsilon_k$ was fitted to data $\{(\hat\sigma_{0,k},\bar F_k):k=1,\ldots,K\}$ 
and yielded the estimates of $\gamma_0$ and $\gamma_1$. 

When applying the same value of tuning parameters $R_k$ in all subwindows, we observed that the estimates of $\gamma_0$ and $\gamma_1$ were rather robust with respect to the choice of $R_k$ in the range from 6 to 11 meters.
In order to include as many pairs of points as possible into the composite likelihood criterion (recall that each subwindow contains only approximately 141 points), we finally decided to use the value $R_k=11$ meters for all subwindows.

 Step 3: For the estimation of $\alpha$, we considered $\sigma(u) = \hat\gamma_0 + \hat\gamma_1 Z(u)$ as fixed and maximized the second-order composite likelihood \eqref{e:CL2} with respect to $\alpha$ only, thereby obtaining an estimate $\hat\alpha$. Note that since $\sigma(u)$ is non-constant, we could not take advantage of the Lebesgue-Stieltjes integral~\eqref{eq:LSint} but needed in each iteration step to calculate the double integral in \eqref{e:CL2}. 
 For the choice of the tuning parameter $R$, 
preliminary experiments showed that the estimates of $\alpha$ were rather stable for $R$ between 5 and 15 meters. Since the width of the observation window $W$ is only 10 meters, it does not seem reasonable to use a very large value of $R$. We decided to use the same value as in the second step, that is, $R=11$ meters.

\subsection{Results}\label{s:resultsfish}

\begin{figure}[t]
\centering 
  \includegraphics[width=\textwidth]{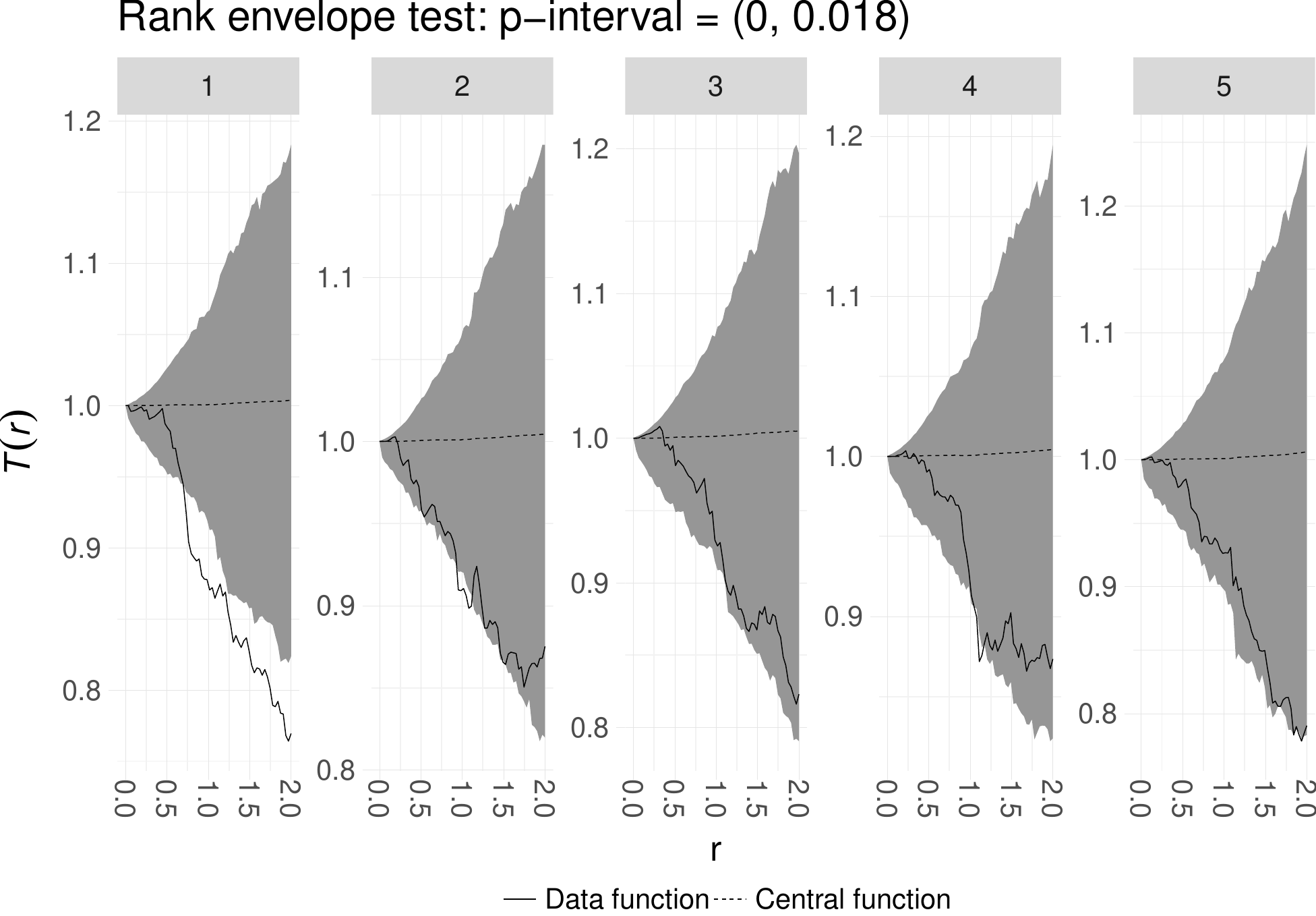} 
   \caption{Rank envelope test for the null hypothesis of no interaction for the fish aggregation dataset.}\label{f:envelope_fish}
\end{figure}

Having estimated the intensity function in the first step of the estimation procedure, it is straightforward to test the null hypothesis of no interaction between the individual fish by performing the rank envelope test and to construct the envelope based on simulations of the fitted inhomogeneous Poisson process. Figure~\ref{f:envelope_fish} shows the simultaneous 95\% envelope constructed from 9999 simulated curves, when the test statistic $T$ is given by concatenating together the five inhomogeneous $J$-functions estimated in the subwindows $W_1, \ldots, W_5$. The data curve (solid black line) leaves the envelope for four of the five subwindows, and the $p$-value is below 2\%, so we reject the null hypothesis of no interaction between the fish. Figure~\ref{f:envelope_fish} also shows that the level of clustering is bigger for shallower subwindows. 

\begin{figure}[t]
\centering 
  \includegraphics[width=9cm]{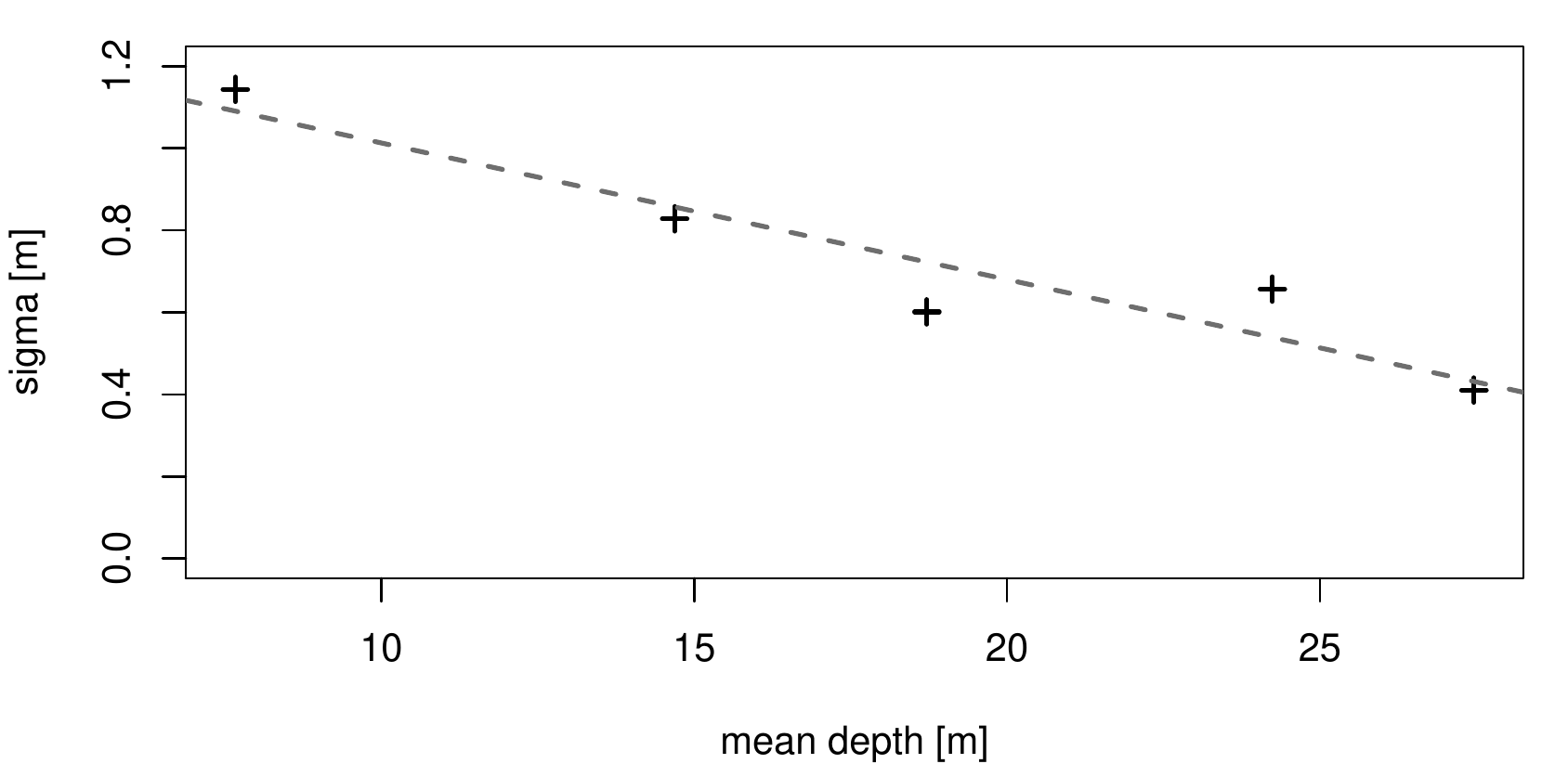} 
   \caption{Estimated values of $\sigma_{0,k}$ plotted against the mean value of the covariate in the subwindow $W_i, i = 1, \ldots, 5$, together with the fitted regression line.}\label{f:regression_fish}
\end{figure}

In the second step of the estimation procedure, a linear regression of $\hat\sigma_{0,1}, \ldots, \hat\sigma_{0,5}$ against average water depths, see Figure~\ref{f:regression_fish}, yields the estimates $\hat\gamma_0$ and $\hat\gamma_1$ of the parameters governing the second-order non-stationarity of the point process. By ignoring uncertainties in the estimates $\hat\sigma_{0,1}, \ldots, \hat\sigma_{0,5}$, the parameter $\gamma_1$ is significantly different from zero ($p$-value=0.015; value obtained with the \texttt{R} function \texttt{glm}).

Table~\ref{tab:fish.res} provides the estimates of the model parameters, including $\hat\alpha$ obtained in the third step. The table also specifies central 90\%-percentile intervals of the estimates obtained with a parametric bootstrap approach \citep{efron1993} where we re-estimated the model parameters from each of 1000 simulated point pattern datasets obtained under the fitted LGCP. Using this parametric bootstrap approach here is more demanding than for the particle-dispersal application considered later, since the double integral in \eqref{e:CL2} needs to be evaluated at each iteration of the maximization procedure required by the 3rd step. A typical computation time took about a minute on a regular desktop computer for datasets with about 700 points. 

\begin{figure}[t]
\begin{center}
\includegraphics[width=\textwidth]{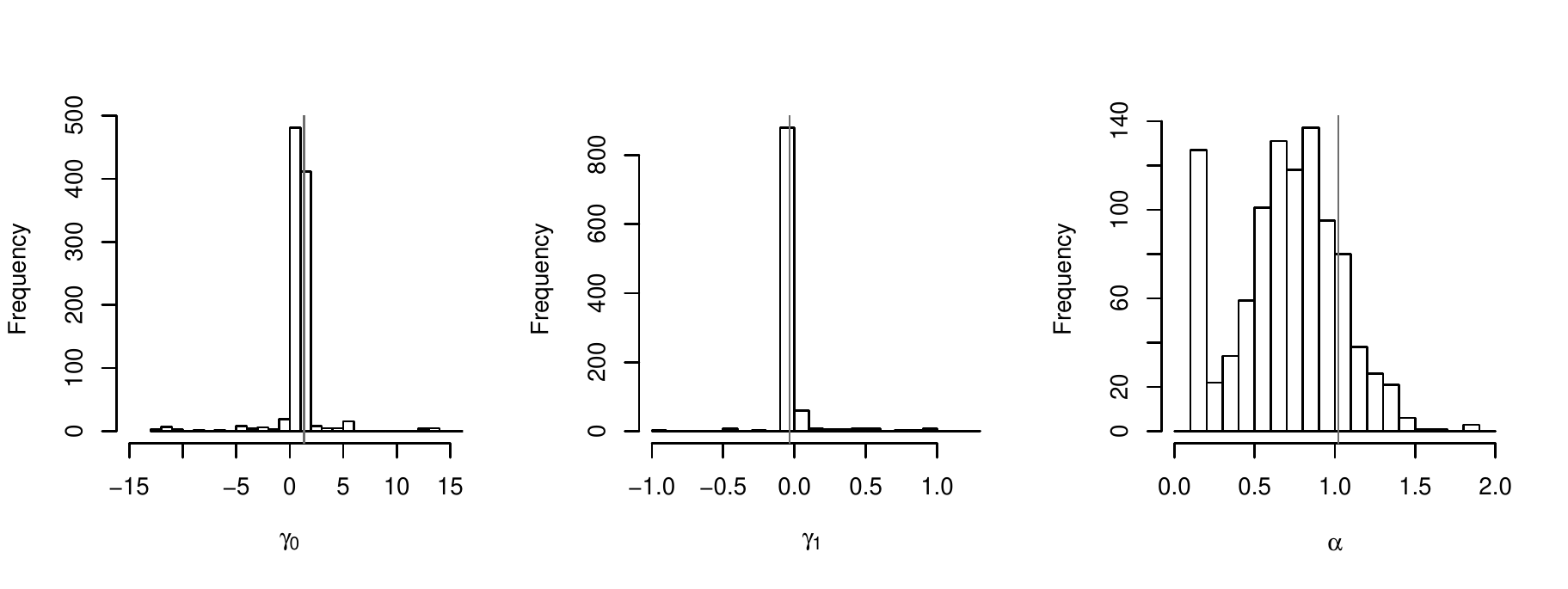} 
\end{center}
\caption{Parameter estimates of the LGCP fitted to the fish aggregation dataset (grey vertical lines) and their bootstrap-based marginal distributions (histograms).}
 \label{fig:fish.estim}
\end{figure}

In addition to Table~\ref{tab:fish.res}, Figure \ref{fig:fish.estim} shows the empirical marginal distributions of the estimates obtained from the 1000 simulated datasets. The regression parameters $\gamma_0$ and $\gamma_1$ governing the variance structure of the GRF $Y_0$ have empirical distributions with a high concentration around their estimated values and a small amount of extremes. The empirical distribution of the scaling parameter $\alpha$ includes a high proportion of values near zero, corresponding to a long-range dependence in the GRF $Q$ and hence also in the LGCP, cf.\ \eqref{e:logL}.

\begin{figure}[t]
\centering 
  \includegraphics[width=\textwidth]{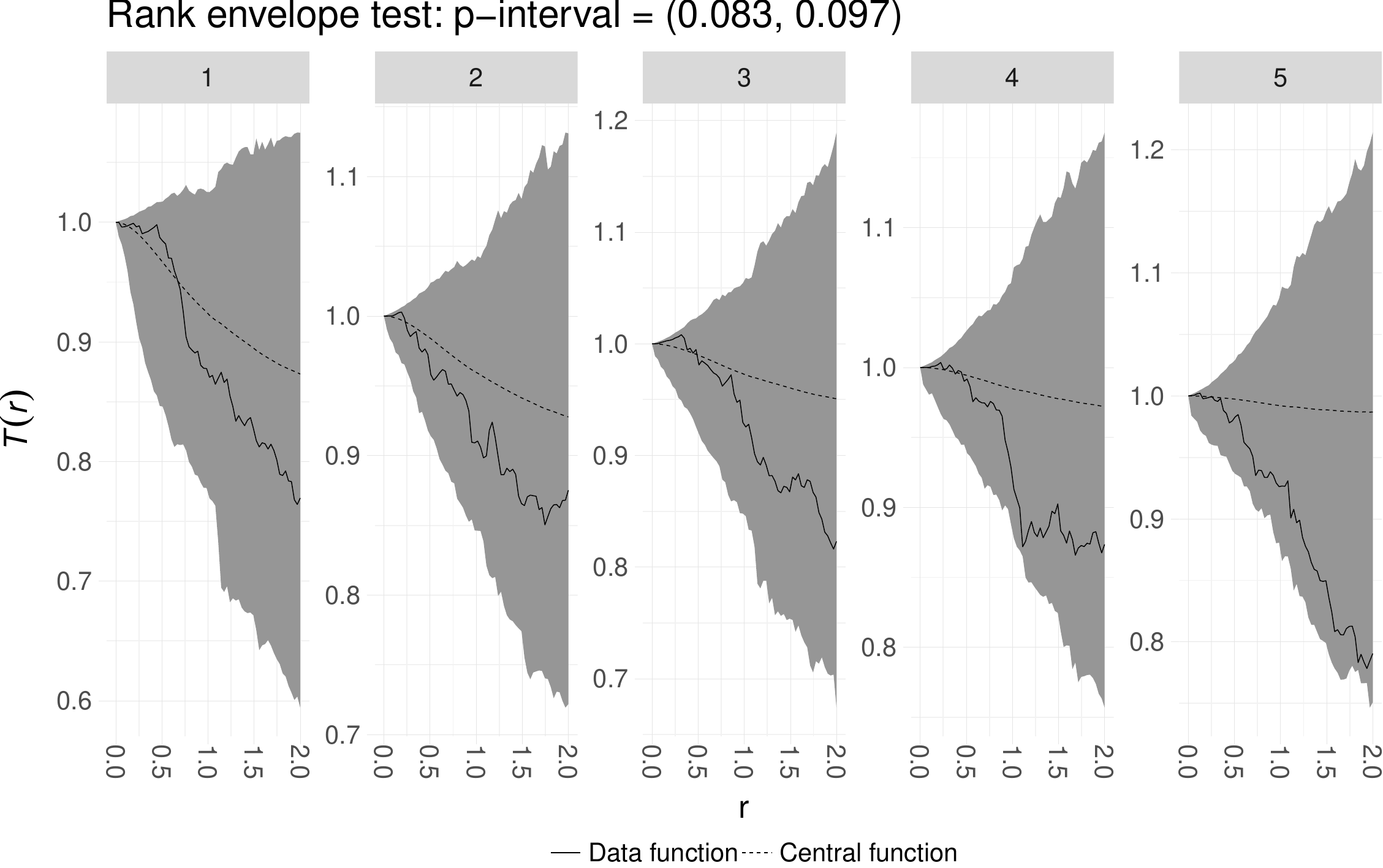} 
   \caption{Rank envelope test to assess the goodness-of-fit of the LGCP fitted to the fish aggregation dataset. Envelopes and the $p$-interval were computed with 9999 simulations.}\label{f:envelope2_fish}
\end{figure}

Goodness-of-fit of the fitted LGCP was assessed with the global envelope test described in Section~\ref{sec:model-check}. For this we used 9999 simulations to obtain the envelopes of the concatenated inhomogeneous $J$-functions computed over the five subwindows $W_1,\ldots,W_5$. The fitted LGCP model is not rejected by the test at the 5\% risk level ($p$-interval=[0.083,0.097]); see Figure~\ref{f:envelope2_fish}. The central (dotted) curve of the global envelope test reveals a decrease in the level of clustering with respect to increasing water depth.

\begin{table}[ht]
\centering
{\small
\begin{tabular}{ccc}
  \hline
  Parameter & Estimate & 90\%-CI \\
  \hline
$\gamma_0$ & 1.35 & [-0.17 , 1.81]\\
$\gamma_1$ & -0.033 & [-0.057 , 0.076]\\
$\alpha$   & 1.02 & [0.10 , 1.22] \\
   \hline
\end{tabular}
}
\caption{\label{tab:fish.res} Estimated parameter values of the LGCP fitted to the fish aggregation dataset, together with their central 90\%-percentile intervals obtained by a parametric bootstrap approach as discussed in Section~\ref{s:Tomasdata}.}
\end{table}

\begin{figure}[t]
\centering 
  \includegraphics[width=0.9\textwidth]{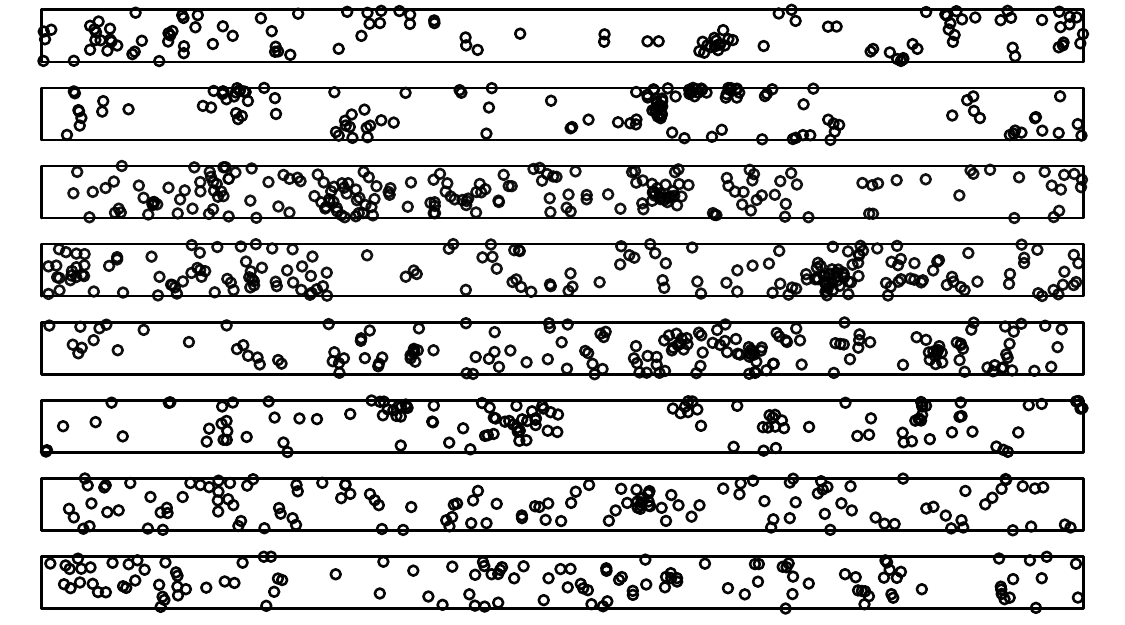} 
   \caption{Simulated realizations from the study in Section~\ref{s:resultsfish}; only a 200 meter long part of the data is shown and the width of the observation window is 10 meters. From top to bottom: Series of simulations 1 to 8, see also Table~\ref{t:JM1} for the corresponding parameter values. The function $\sigma(u)$ is constant when $\gamma_1 = 0$ (Series 1, 2, 5, and 6) and  varying when $\gamma_1 = -0.03$ (Series 3, 4, 7, and 8); the range of spatial autocorrelation is smaller for $\alpha=0.50$ (Series 5--8) and larger for $\alpha=0.25$ (Series 1--4); the variance of the underlying random field is smaller for $\gamma_0=1.25$ (Series 1, 3, 5, 7) and larger for $\gamma_0=1.50$ (Series 2, 4, 6, 8).}\label{f:simulated_fish}
\end{figure}

\begin{figure}[t]
\begin{center}
\includegraphics[width=\textwidth]{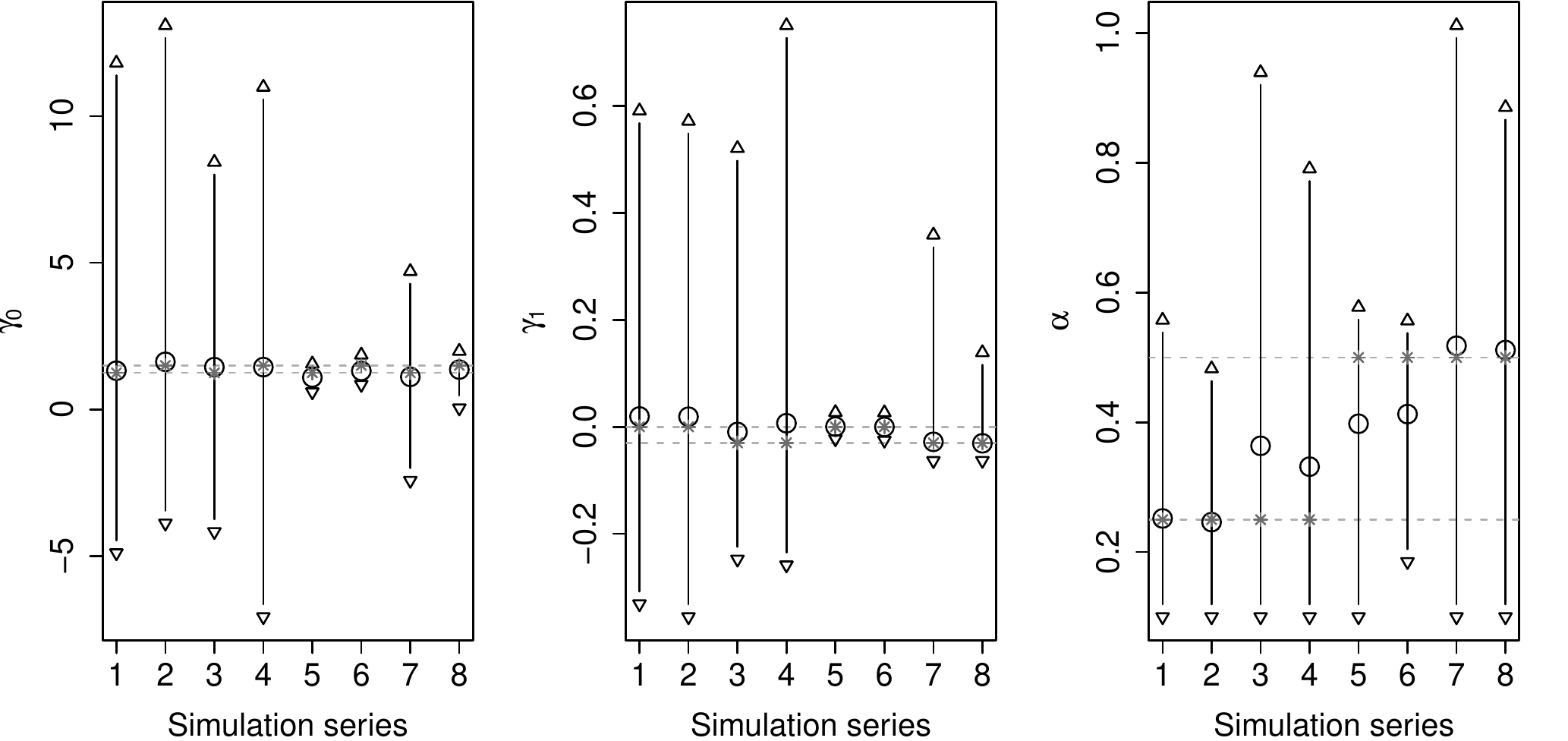} 
\end{center}
\caption{Summaries of the empirical distributions of parameter estimates for the LGCP with the variance of the random field depending on a covariate and fitted to 8000 simulated datasets which are split into 8 series of simulations corresponding to different parameter values (see main text). Stars: true parameter values. Circles: trimmed means of parameters obtained by removing the lower and upper 5\% of estimates (the trimming was done to moderate the effect of extreme estimates). Triangles: endpoints of 90\%-confidence intervals.}
 \label{fig:fish.simul}
\end{figure}

We investigated the efficiency of our estimation procedure for 8 series of 1000 simulations corresponding to the 8 combinations of parameter values obtained when $\alpha$, $\gamma_0$, and $\gamma_1$ each takes two values, see Table~\ref{t:JM1}. Figure~\ref{f:simulated_fish} shows examples of the simulated realizations. Note that the function $\sigma(u)$ is constant when $\gamma_1 = 0$ (Series 1, 2, 5, and 6) and  varying when $\gamma_1 = -0.03$ (Series 3, 4, 7, and 8), and our choice of parameter values was in part constrained by the range of the covariate values (5 to 29 meters) and made to avoid negative values for $\sigma(u)$ or values close to zero for $\sigma(0)$. 
Further, in the simulations, we used the same observation window $W$ and the same covariate describing the water depth as for the fish aggregation dataset, and also we applied the same estimation procedure. In particular, we used $R_k=11$ for all $k$.
Figure~\ref{fig:fish.simul} provides estimation results for each series of simulations. The regression parameters $\gamma_0$ and $\gamma_1$ governing the spatially structured variance function of the random field are estimated with negligible bias whatever the value of $\alpha$, and with rather small variability when the spatial correlation of the random field is small ($\alpha=0.50$; Series 5--8). On the other hand the estimation of $\alpha$ is more difficult. This is not surprising because $\alpha$ is estimated in the third estimation step and relies on the estimated values of $\gamma_0$ and $\gamma_1$.

\begin{table}[t]
\centering
{\small
\begin{tabular}{ccccccccc}
  \hline
  
  & \multicolumn{8}{c}{Series of simulations}\\
 & 1 & 2 & 3 & 4 & 5 & 6 & 7 & 8 \\
  \hline
$\alpha$ & 0.25 & 0.25 & 0.25 & 0.25 & 0.50 & 0.50 & 0.50 & 0.50\\
$\gamma_0$ & 1.25 & 1.50 & 1.25 & 1.50 & 1.25 & 1.50 & 1.25 & 1.50\\
$\gamma_1$ & 0 & 0 & -0.03 & -0.03 & 0 & 0 & -0.03 & -0.03 \\
   \hline
\end{tabular}
}
\caption{\label{t:JM1} Parameter values used in the simulation study related to the estimation procedure for the fish aggregation dataset.}
\end{table}

Table \ref{T_fish_test} shows the results of the test of the hypothesis $\gamma_1=0$ based on the linear regression of $\hat\sigma_{0,1}, \ldots, \hat\sigma_{0,5}$ against the covariate used to model $\sigma$ (namely, the average water depths in the application). This test, which is a way to assess the dependence of the second order structure of the LGCP model on the specified covariate, was applied to the 8 series of 1000 simulations with different levels of significance (1\%, 5\% and 10\%). Based on the results for Series 1, 2, 5, and 6, the actual significance level of the test is approximately equal to or lower than the nominal significance level, that is, the test is slightly conservative. The power of the test, assessed with Series 3, 4, 7, and 8, is increased when the spatial auto-correlation is larger at shorter distances (smaller $\alpha$), but could certainly be improved with a test better accounting for the propagation of uncertainties in the first two steps of the estimation procedure.

\begin{table}[ht]
\centering
{\small
\begin{tabular}{ccccccccc}
  \hline
  
Significance  & \multicolumn{8}{c}{Series of simulations}\\
level & 1 & 2 & 3 & 4 & 5 & 6 & 7 & 8 \\
  \hline
 0.01& 6 &10& 56& 57& 4& 5& 42& 40\\
0.05& 40& 50& 236& 254& 33& 23& 196& 156\\
0.10& 77& 101& 359& 410& 62& 64& 305& 284 \\
   \hline
\end{tabular}
}
\caption{\label{T_fish_test} Numbers of rejections (over 1000) in the test of $\gamma_1=0$ applied to the 8 series of simulations with different significance levels.}
\end{table}

\section{Particle dispersal}\label{sec:particles}

\subsection{Data and model}\label{s:Samueldata}

In this section, we are interested in point patterns generated by the deposit locations of biological windborne particles (such as seeds, pollen grains and pathogenic fungal spores) dispersed from plants. In general, for an isolated plant releasing particles, the point pattern formed by the deposit locations of the particles is denser around the source plant than far from it. In other words, the intensity of points generally decreases with the distance from the source along radial directions \citep{austerlitz2004,rieux2014,soubeyrand2007phyto,tufto1997}. This large-scale inhomogeneity may be augmented by small-scale variations due, for instance, to group dispersal \citep{soubeyrand2011}: groups of particles are transported independently, but the particles of a given group are deposited at dependent locations and form a cluster of points whose spatial extent is smaller at shorter distances (dependencies between the particles of a group vanishes at longer distance because of the accumulation of air turbulence; see~\ref{app:dispersal.motiv} for an extended explanation). In such a case, deposit locations of particles form an inhomogeneous point pattern with spatially varying local aggregation, which may be modelled by a LGCP with a non-stationary auto-correlation function for the GRF.
We used this modelling framework and the accompanying estimation procedure to the point pattern in Figure~\ref{f:particles} formed by the deposit locations of fungal spores dispersed from a point source located at the origin \citep[details about this data set can be found in][]{lannou2008}. Here, we assume a log linear intensity as in \eqref{e:log-linear-rho} with $I=1$, $z_1(u)=\|u\|$, and $(\beta_0,\beta_1)\in\mathbb R^2$. Such a specification corresponds to the so-called  exponential dispersal kernel $u\rightarrow\exp(\beta_0)\exp(\beta_1\|u\|)$, and it means that the intensity of deposit locations of particles decreases exponentially with distance from the point source when $\beta_1<0$. Let $\sigma(u)\equiv \sigma$ be a constant function, and use the space transformation 
\begin{equation}\label{eq:phi}
\varphi(u)=\log(1+\delta\|u\|)\|u\|^{-1}u,
\end{equation}
where $\delta$ is a non-negative parameter for modelling an eventual increase in the spatial extent of clusters with the distance from the source point. Finally, assume an exponential correlation function of the GRF:
\begin{equation}\label{e:expcorr}
r(t)=\exp\left(- \alpha t\right),\qquad t\ge0,
\end{equation}
where $\alpha$ is a positive parameter.

\subsection{Implementation of the three step procedure}\label{s:parest-samuel}

We adapted the three step estimation procedure of Section \ref{sec:3-steps} to the LGCP described above incorporating the parametric space transformation \eqref{eq:phi}. Here we let $\nu^\text{split}=\delta$ and $\nu^\text{joint}=(\sigma,\alpha)$.

Step 1: For the estimation of $\beta_0$ and $\beta_1$, we maximized the first-order composite likelihood to obtain estimates $\hat\beta_0$ and $\hat\beta_1$. 

Step 2: For the estimation of $\delta$, we started by defining subwindows $W_1,\ldots,W_K$ based on the distance from the point source by setting $f(u)=\|u\|$ so that $W_1$ is the inner disk and $W_2,\ldots,W_K$ are the subsequent rings; more details are given below. 
Then, for each $k=1,\ldots,K$, we estimated $\nu_{0,k}=(\sigma_{0,k},\alpha_{0,k})$ by  maximizing the approximate second-order log composite likelihoods based on $X\cap W_k$. Let $\hat\alpha_{0,k}$ denote the estimate of $\alpha_{0,k}$.  Assuming that the function $f$ is approximately equal to a constant $\bar F_k$ on each subwindow $W_k$, the covariance function of $Y_0$ can be approximated for $u,v\in W_k$ by
\begin{align*}
c(u,v)&=\sigma^2 \ r(\|\varphi(u)-\varphi(v)\|) \\
&= \sigma^2 \ r\left(\|\log(1+\delta\|u\|)\|u\|^{-1}u - \log(1+\delta\|v\|)\|v\|^{-1}v\|\right)\\
&\approx \sigma^2 \ r\left(\|\log(1+\delta\bar F_k)\bar F_k^{-1}u - \log(1+\delta\bar F_k)\bar F_k^{-1}v\|\right)\\
&= \sigma^2 \ \exp\left(-\alpha \bar F_k^{-1} \log(1+\delta\bar F_k)   \|u-v\|\right).
\end{align*}
Hence, we first estimated $\alpha \bar F_k^{-1} \log(1+\delta\bar F_k) $ by $\hat\alpha_{0,k}$ and second, to obtain an estimate $\hat\delta$ of $\delta$, 
we fitted the non-linear regression model $\alpha_{0,k}=\alpha \bar F_k^{-1} \log(1+\delta\bar F_k)+\epsilon_k$ to the data $\{(\hat\alpha_{0,k},\bar F_k):k=1,\ldots,K\}$ by minimizing $\sum_{k=1}^K | \hat\alpha_{0,k} - \alpha \bar F_k^{-1} \log(1+\delta\bar F_k) |$ with respect to $\alpha$ and $\delta$. For our dataset, we chose $K=5$; since only 5 points then are entering the regression, we have used absolute difference in order to increase robustness of the regression estimate. 
We applied the following rules of thumb for selecting $f_k$ and $R_k$: $f_0,\ldots,f_K$ were selected such that $X\cap W_1,\ldots,X\cap W_K$ contain approximately the same numbers of points and $W=W_1\cup\cdots\cup W_K$; for $k=1,\ldots,K$, we selected $R_k$ such that about 50\% of the observed pairs of points in $W_k$ were used in the computation of the partial approximate second-order log composite likelihood \eqref{e:CL2.split}.

\begin{figure}[t]
\begin{center}
\includegraphics[height=\textwidth,angle=-90]{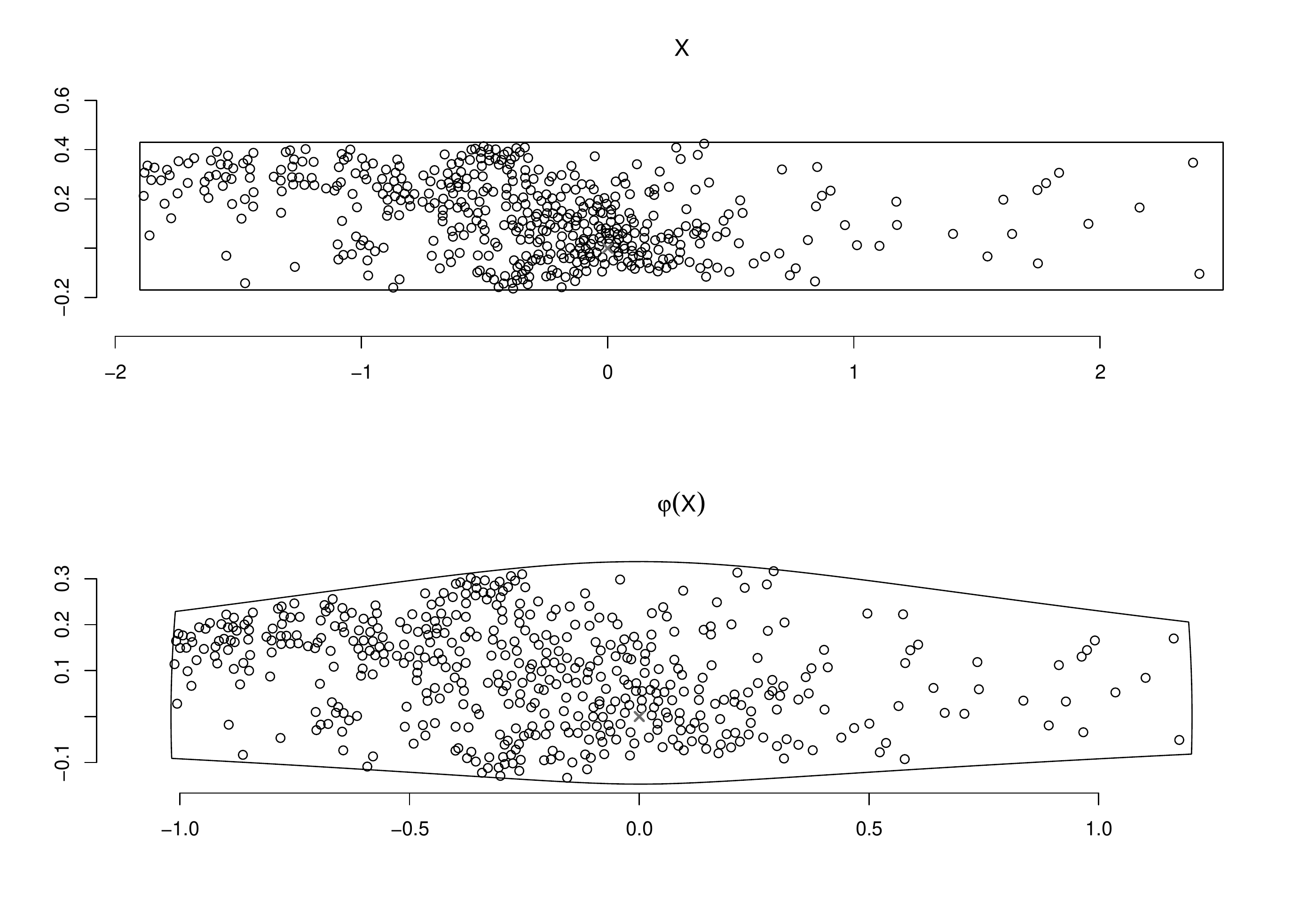}
\end{center}
\caption{Top: Point pattern $X$ (open circles) analysed in the particle dispersal application and location of the source of dispersed particles (cross). Bottom: transformed point pattern $\hat\varphi(X)$.} 
 \label{fig:app.plot.varphiX}
\end{figure}

 Step 3: For the estimation of $\sigma$ and $\alpha$ we used the transformed process $\hat{\varphi}(X)$ where
 $\hat{\varphi}(u) = \|u\|^{-1}\log(1+\hat\delta\|u\|)u$.
 The motivation was that, assuming the true value of $\delta$ is known, $\varphi(X)$ is a second-order intensity-reweighted stationary LGCP \citep{baddeley:moeller:waagepetersen:00} with known form of the intensity function and the pair-correlation function. See \ref{app:varphiX} for a detailed derivation and Figure~\ref{fig:app.plot.varphiX}
 for the original pattern and the pattern transformed by
$\hat{\varphi}(u) = \|u\|^{-1}\log(1+\hat\delta\|u\|)u$,
 where $\hat\delta$ is given in Table~\ref{tab:particles.res}. We approximated the intensity function of $\varphi(X)$ in the observation window $\varphi(W)$ by a log-linear function of $\| u \|$ and estimated it by the first-order composite likelihood. We used this estimated intensity function to construct the second-order composite likelihood criterion for estimation of $\sigma$ and $\alpha$. Again we chose the value of $R$ so that about 50 \% of the observed pairs of points were used in the criterion.

\subsection{Results}\label{sec:disp.res}

\begin{figure}[t]
\begin{center}
\includegraphics[height=\textwidth,angle=-90]{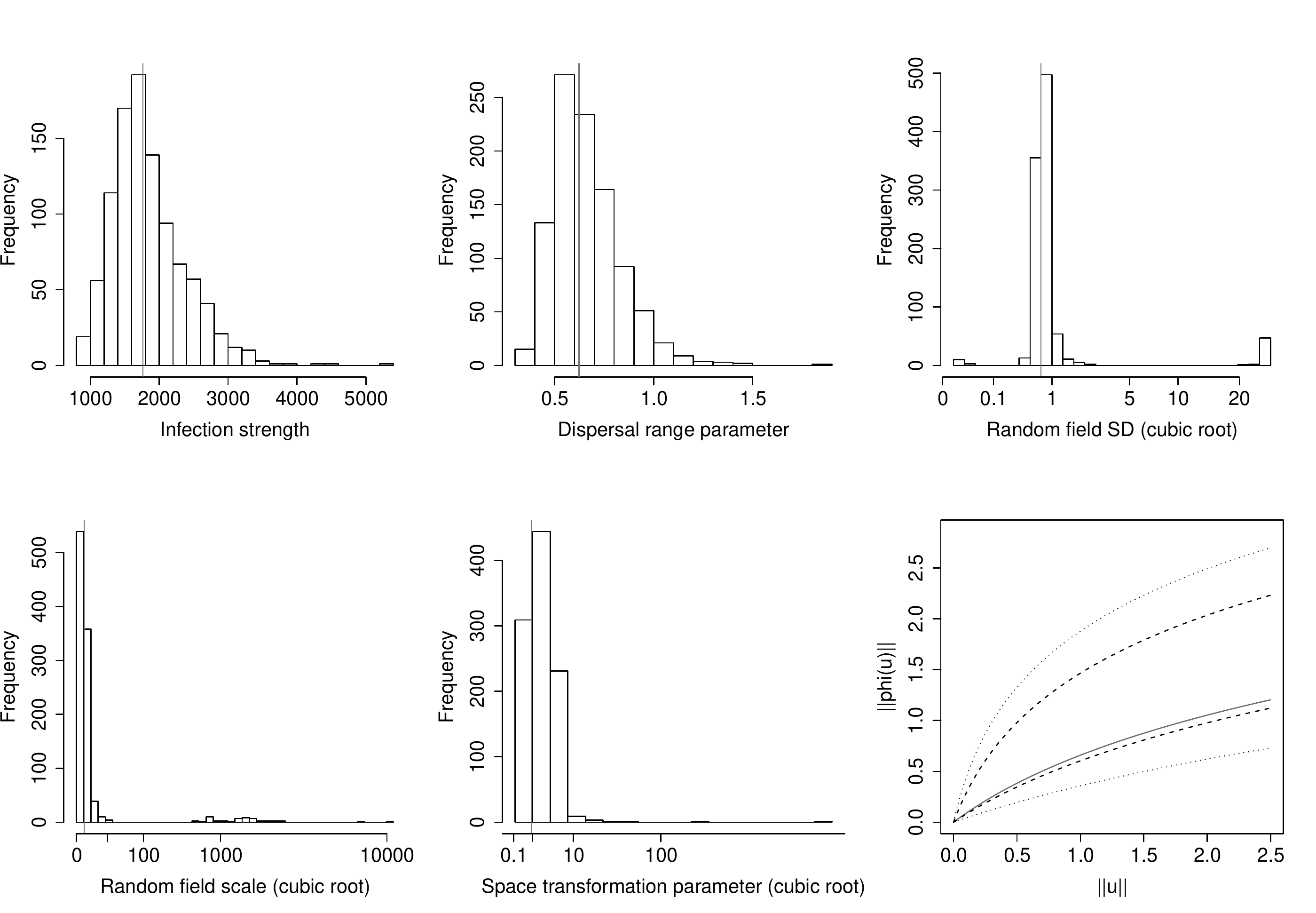} 
\end{center}
\caption{Estimates of transformed parameters of the LGCP with space transformation fitted to the particle dispersal dataset (vertical lines) and their bootstrap-based marginal distributions (histograms); expressions of transformed parameters are provided in Table \ref{tab:particles.res}. The bottom right panel provides the estimated space transformation function $\hat\varphi$ (solid curve) 
 and its pointwise 90\%- and 50\%-confidence envelopes (dotted and dashed lines, respectively).  }
 \label{fig:particles.estim}
\end{figure}

\begin{figure}[t]
\begin{center}
\includegraphics[height=\textwidth,angle=-90]{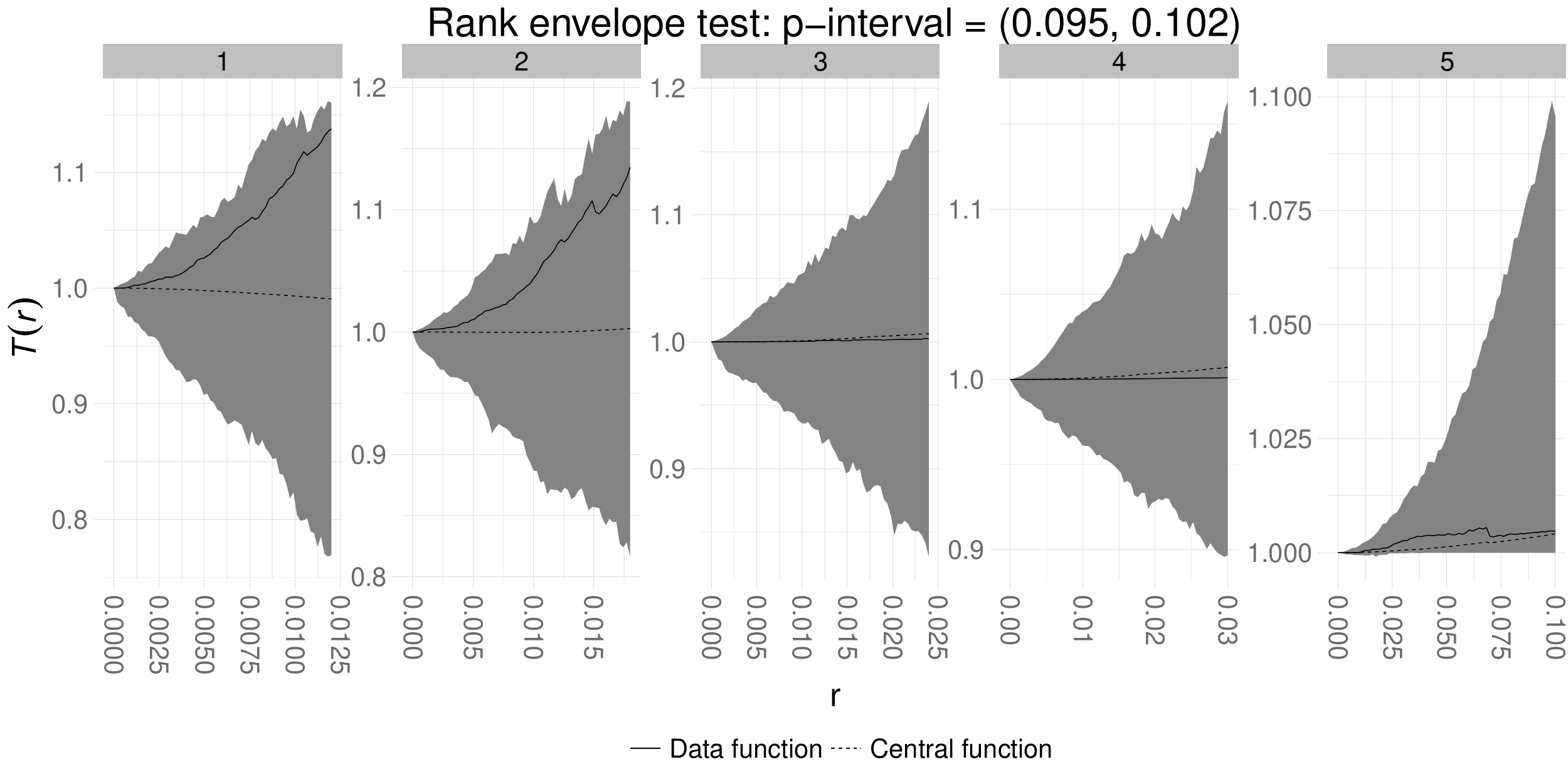}
\end{center}
\caption{Rank envelope test to assess the goodness-of-fit of the LGCP with space deformation fitted to the particle dispersal dataset. Envelopes and the $p$-interval were computed with 9999 simulations.}
 \label{fig:particles.test}
\end{figure}

Table \ref{tab:particles.res} provides the estimates of (transformed) model parameters obtained with the three step procedure. It also gives their 90\%-percentile intervals obtained with the same parametric bootstrap procedure as in Section~\ref{s:resultsfish}; this was fast:
it took typically a few seconds on a regular desktop computer for datasets with about 500 points. Figure~\ref{fig:particles.estim} shows the marginal distributions of estimates obtained from the 1000 simulated datasets, and it shows the estimated space transformation function $\hat\varphi$ and its pointwise 90\%- and 50\%-confidence envelopes.
We notice a rather large estimation variability for second-order parameters $(\sigma,\alpha,\delta)$.
The goodness-of-fit of the model was assessed with the global envelope test described in Section~\ref{sec:model-check}, using 9999 simulations to obtain the envelopes of the  concatenated inhomogeneous $J$-functions computed over the five subwindows $W_1,\ldots,W_5$.
The fitted LGCP model is not rejected by the test at the 5\% risk level ($p$-interval=[0.095,0.102]); see Figure~\ref{fig:particles.test}.
The observed curve indicates repulsion in the first two subwindows. This is probably a consequence of the non-negligible size of the particles, an assumption that is not accounted for in our model.

\begin{table}[ht]
\centering
{\small
\begin{tabular}{lcccc}
  \hline
  Name of transformed parameter & Parameter expression & Estimate & 90\%-CI \\
  \hline
Infection strength & $2\pi e^{\beta_0}/\beta_1^2$ & 1760 & [1140 , 2820]\\
Dispersal range parameter & $-1/\beta_1$ & 0.62 & [0.44 , 0.98]\\
Random field standard-deviation & $\sigma$ & 0.73 & [0.55 , 3.50] \\
Random field scale & $1/\alpha$ & 0.15 & [0.02 , 34.3] \\
Space transformation parameter & $\delta$  & 0.93 & [0.42 , 5.54]\\
   \hline
\end{tabular}
}
\caption{\label{tab:particles.res} Estimation of parameters of the LGCP with space transformation fitted to the particle dispersal dataset. Original parameters are transformed to obtain interpretable quantities. The infection strength gives the expected total number of spores released by the source. The dispersal range parameter is half the mean dispersal distance of spores. The random field scale parameter is the distance at which the auto-correlation is $1/e$ when there is no space transformation.}
\end{table}

\begin{figure}[t]
\centering 
  \includegraphics[angle=90,width=0.9\textwidth]{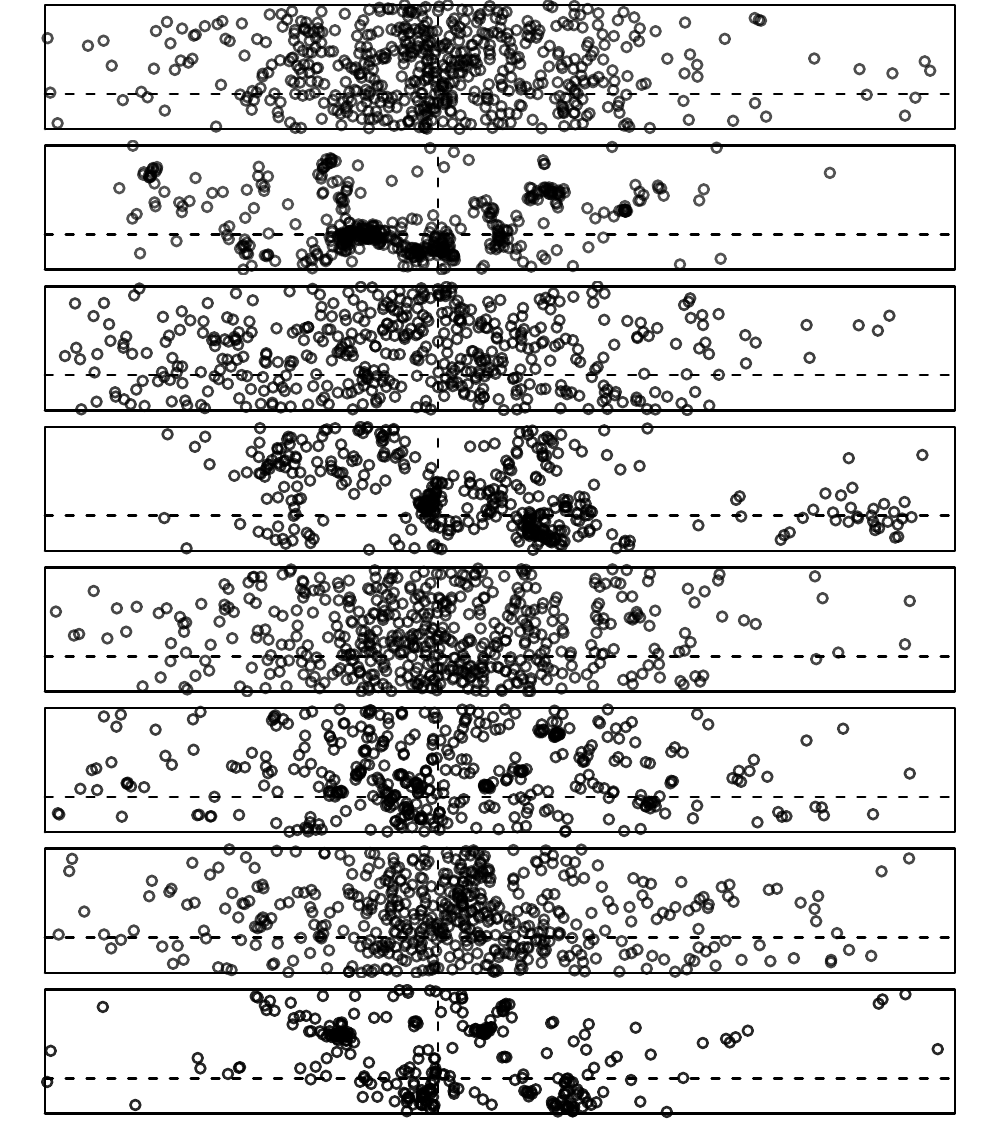}
   \caption{Simulated realizations from the study in Section~\ref{sec:particles}. From left to right: Series of simulations 1 to 8, see also Table~\ref{t:JM2} for the corresponding parameter values. Dashed lines represent the coordinate axes. The variability of the underlying random field is smaller for $\sigma = 0.5$ (Series 1, 3, 5, 7) and larger for $\sigma = 2$ (Series 2, 4, 6, 8); the range of spatial autocorrelation is smaller for $1/\alpha = 0.05$ (Series 1, 2, 5, 6) and larger for $1/\alpha=0.2$ (Series 3, 4, 7, 8); the space deformation is very prominent for $\delta = 10$ (Series 5--8) and not so prominent for $\delta=1$ (Series 1--4).}\label{f:simulated_GD}
\end{figure}

\begin{figure}[t]
\begin{center}
\includegraphics[height=\textwidth,angle=-90]{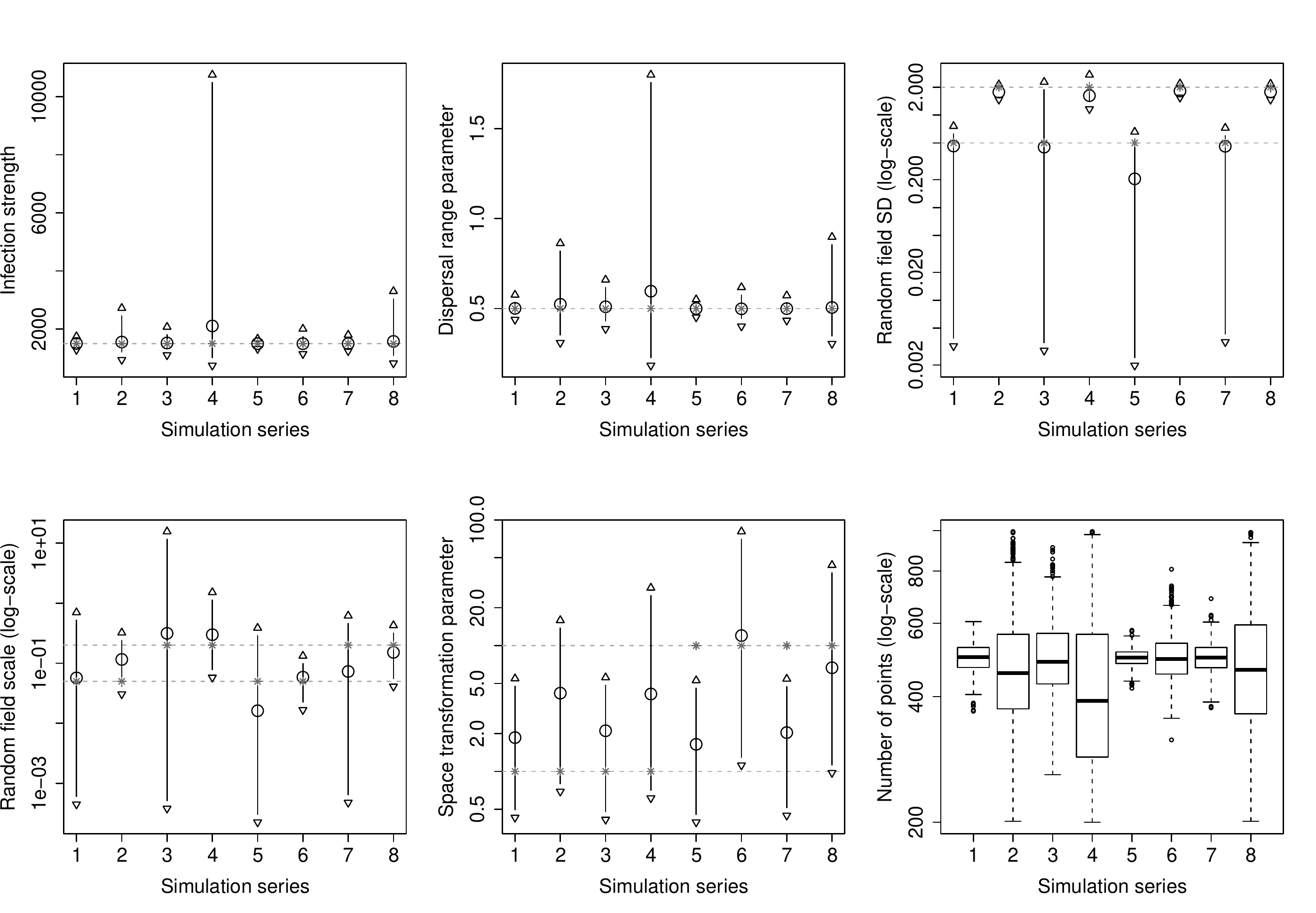}
\end{center}
\caption{Estimation of transformed parameters of the LGCP with space transformation fitted to 8000 simulated datasets split into 8 series of simulations corresponding to different parameter values (see main text); expressions of transformed parameters are provided in Table \ref{tab:particles.res}. Stars: true parameter values. Circles: trimmed means of parameters obtained by removing the lower and upper 5\% of estimates (the trimming was done to moderate the effect of extreme estimates). Triangles: endpoints of 90\%-confidence intervals. The bottom right panel provides the distribution of the number of points in $W$ for each simulation series. This number of points was constrained to be between 200 and 1000.}
 \label{fig:particles.simul}
\end{figure}

We investigated the efficiency of our estimation procedure for 8 series of 1000 simulations performed with 8 different sets of parameter values. In all series, infection strength and dispersal range were equal to $(2\pi e^{\beta_0}/\beta_1^2,-1/\beta_1)=(1500,0.5)$, whilst the values of $(\sigma,1/\alpha,\delta)$ are given in Table~\ref{t:JM2}. Figure~\ref{f:simulated_GD} shows examples of the simulated realizations.
Figure~\ref{fig:particles.simul} provides estimation results and the distribution of the number of points for each simulation series. 
The estimation of the infection strength and the dispersal range is quite accurate except when the GRF $Y_0$ is strongly variable ($\sigma=2$) and its spatial correlation is large ($1/\alpha=0.2$). The latter situation of the GRF can generate large clusters of points in the LGCP away from the source of particles and negatively impact the estimation accuracy of the spatial trend in the intensity of points, which is supposed to be decreasing with the distance from the source of particles. Moreover, if the GRF's standard-deviation is appropriately estimated, we observe rather large estimation variability and possible biases for $1/\alpha$ and $\delta$, the parameters arising in the non-stationary auto-correlation function. Overall, the inference is relatively accurate when the GRF's variance is large ($\sigma=2$) and the spatial deformation is significant ($\delta=10$); i.e.,  series 6 and 8. In contrast, for a rather flat GRF or when the deformation is not so significant, the estimation of parameters governing the second-order structure is clearly uncertain due to a complex model and the limited amount of information in the data.

In \ref{app:comparison} we compare the performance of the three step estimation method described above with a two step procedure using in the second step the composite likelihood criterion based on the whole observation window $W$. We conclude that the three-step procedure is in all cases more successful in estimation of $\sigma$ and in almost all cases for $\alpha$. For the deformation parameter $\delta$ the comparison is less clear but in most cases the three-step procedure results in more precise estimates.

\begin{table}[ht]
\centering
{\small
\begin{tabular}{ccccccccc}
  \hline
  
  & \multicolumn{8}{c}{Series of simulations}\\
 & 1 & 2 & 3 & 4 & 5 & 6 & 7 & 8 \\
  \hline
$\sigma$ & 0.5 & 2 & 0.5 & 2 & 0.5 & 2 & 0.5 & 2\\
$1/\alpha$ & 0.05 & 0.05 & 0.2 & 0.2 & 0.05 & 0.05 & 0.2 & 0.2\\
$\delta$ & 1 & 1 & 1 & 1 & 10 & 10 & 10 & 10 \\
   \hline
\end{tabular}
}
\caption{\label{t:JM2} Parameter values used in the simulation study related to the estimation procedure for the particle dispersal dataset.}
\end{table}

\section{Concluding remarks}\label{s:conclusion}

In this article, we modelled inhomogeneity in spatial point patterns by a Cox process driven by a second-order non-stationary log GRF. We presented a first example where the non-stationarity was obtained by modelling the variance of the GRF as a linear function of an explanatory variable. In a second example, the non-stationarity was obtained via a spatial parametric deformation,  
which enables the generation of point clusters with spatially varying spatial extents. In this framework, estimating model parameters, including those that govern the second-order non-stationarity, is challenging, especially with a moderate number of points (a few hundreds in our examples). We developed and tested a three step estimation approach based on first- and second-order composite likelihoods.
The second-order composite likelihood has to be specifically adapted to the non-stationarity incorporated into the model. The advantages of this estimation procedure are the reduced computation time and the reduced implementation complexity. Reduced computation time allowed us to implement parametric bootstrap 
for assessing estimation uncertainty, and computationally intensive model check approaches,
namely global envelope tests,
for assessing goodness-of-fit.

We leave it for future research to examine the theoretical properties of the estimators in the three step estimation procedure. However, the performance of this procedure was evaluated in two simulation studies with settings inspired by the two case studies. Our simulation studies highlight situations where a satisfactory estimation accuracy can be obtained and, conversely, situations where estimation performance is poor. Poor estimation mostly translates into large estimation variance and, from a general point of view, arises when variations in the underlying random field are weak or the spatial deformation is not so significant (due to a complex model and the limited information contained in the data).
A part of the estimation poorness can also be attributed to the choice of a three steps procedure, which neglects uncertainty propagation. However, this procedure definitely allows to reduce computation time. Furthermore, it has to be noted that tuning parameters used in the calculation of the second-order composite likelihood have to be appropriately specified to reach satisfactory estimation performance (see Sections \ref{s:parest-tomas} and \ref{s:parest-samuel}). 
Generally, one can proceed as follows. Once the model is chosen, the number of subwindows has to be decided. This can be done in order to have in each subwindow sufficiently many points for reliable estimation with the two step procedure for the second-order intensity-reweighted stationary case (on the union of subwindows) and on the other hand to have subwindows small enough to be able to assume second-order intensity-reweighted stationarity within each subwindow. About 100 points in each subwindow could be recommended as a first choice but some experimenting with this choice is wise.  Then the maximum distances $R_k$ has to be chosen. It is possible to use the same $R_k$ values for all subwindows in cases where the interaction distance is homogeneous, cf.\ Section~\ref{sec:fish}. It is also possible to use $R_k$ values such that they take into account a given percentage of the pairs of points in cases where the interaction distance is inhomogeneous, cf.\ Section~\ref{sec:particles}. Then, in both cases, we recommend to choose the maximum $R_k$ which gives stable results, i.e., similar to those obtained with smaller $R_k$. Finally, for the third step we recommend to use the same maximum distance.

In order to graphically explore the dependence of the second order structure on a given covariate we utilized the global envelope test of independence of points locations, which was performed on first order inhomogeneous $J$-function computed across several subwindows which respect to the value of the given covariate (e.g.\ Figure \ref{f:envelope_fish}). 

We also introduced an approximate test of dependence of the second order structure on a given covariate based on regression of second order parameters estimated across the subwindows which respect to the value of the given covariate. Our simulation study shows that this test is slightly conservative, which is acceptable in practise.

As discussed in the first paragraph of this section, the two studied examples give two different approaches for modelling complex point process structures. But the way of modelling the second-order non-stationarity is not the only difference. Indeed, the first order intensity is described by a nonparametric model in the fish application, and a parametric model in the particle dispersal application. This choice is related to a generic difficulty in modelling, which has to be carefully investigated for every particular dataset: how to divide data variability between the first order intensity and the second order structure. In the case of a nonparametric first order intensity, this choice is done in the selection of the bandwidth. In the case of a parametric first order intensity, this choice is constrained by the flexibility of the parametric model. For instance, using a small bandwidth in the former case can cause overfitting of the first order intensity with no detection of second order structure; and using a too simple parametric model of the intensity in the latter case can induce a bias, which then impacts the estimation of the second order structure. Such issues generally rely on the usual modelling strategy of the analyst but should also rely on the classical process of statistical analysis consisting of defining a class of models, selecting a model based on the data (including parameter estimation), checking model fit, and possibly re-defining the class of models for starting again the process and finally drawing conclusions \citep[Chapter 12]{mccullagh1989}.

Finally, from a prospective point of view, with increasing amount of the gathered data, there is increasing demand for the more complex modelling. The next level for complexity of the presented models is certainly the shift into spatio-temporal modelling framework of e.g. generations of the infection spread or several spatial snapshots of the point pattern gathered in different times. We also leave this issue for future research, noticing that a space time observation window $W\times T$ can be divided in different ways, using sets $W \times T_k, W_k \times T$ or $W_k \times T_j$, depending on the properties of the particular dataset.

\appendix

\section{Small-scale variations in particle dispersal} \label{app:dispersal.motiv}

Section \ref{sec:particles} deals with point patterns generated by the deposit locations of biological windborne particles dispersed from plants. As indicated in the main text, the intensity of points generally decreases with the distance from the source along radial directions. Nevertheless, in particular cases,  the intensity of points along radial directions may be non-monotonic in the vicinity of the source \citep{stoyan2001}. This leads to a large-scale inhomogeneity in the point pattern formed by the deposit locations of the particles. In addition, the inhomogeneity of the point process may be augmented by small-scale variations due, for instance, to dependencies in the dispersal of particles \citep{soubeyrand2011} or to a spatial heterogeneity of the environment \citep{soubeyrand2007jds}. 

Consider a situation where aggregation is due to dependencies in the dispersal of particles, cf.\ Section \ref{sec:particles}. \citet{soubeyrand2011} investigated such a situation, in which (i) the particles are released as groups of particles, (ii) the particles of each group are transported together by the wind but the accumulation of air turbulences progressively separate the particles, and (iii) the particles of each group are deposited at nearby locations but the proximity of the deposit locations depends on the accumulation of air turbulences. Typically, the proximity of the deposit locations decreases with the distance separating the source and the final group center location (the larger this distance is, the larger the accumulation of air turbulences and larger the separation of the particles is). A more complete bio-physical justification for such a process is provided in \citet{soubeyrand2011}. In this situation, one expects that the deposit locations of the particles of a group form a cluster whose spatial extent increases with the distance from the source. Consequently, aggregates of points at longer distances from the source should have wider spatial extent. To model this, we can consider a non-stationary auto-correlation function for the zero-mean and unit-variance GRF $Q$ in \eqref{e:logL}.

In contrast, consider a situation where aggregation is due to a spatial heterogeneity of the environment. If we assume in addition that the spatial heterogeneity is stationary, then the probabilistic properties of aggregates should be the same across space and the auto-correlation function for $Q$  may be stationary.

In the application, the bootstrap analysis is only moderately in favour of the presence of deformation (see bottom right panel of Figure \ref{fig:particles.estim}, where the 90\%-confidence envelope is quite wide). Thus, we cannot conclude that `group dispersal occurs' or that `the spatial extent of groups varies with distance at the scale of the study domain if group dispersal occurs'. Replicated data could be helpful for more deeply investigating this issue.

\section{Properties of $\varphi(X)$}\label{app:varphiX}

Let the situation be as in Section~\ref{sec:particles}. This appendix shows that $\varphi(X)$ is a second-order intensity-reweighted stationary LGCP.

Let $\Lambda(u), u \in \mathbb{R}^2,$ denote the random intensity function of the process $X$ given by $\log \Lambda (u) = \zeta(u) + Y_0 (u), u \in \mathbb{R}^2$, and let $M$ be the corresponding random intensity measure, that is, $M(A) = \int_A \Lambda (u) \, \mathrm{d}u$ for any Borel set $A \subset \mathbb{R}^2$. 
In the third step of the estimation procedure for the particle dispersal dataset, we consider the transformed process $\varphi(X)$. To find the distribution of $\varphi(X)$, we compute void probabilities: For any Borel set $A \subset \mathbb{R}^2$,
\begin{align*}
  \mathrm{P} ( \varphi(X) \cap A = \emptyset ) & = \mathrm{P} ( X \cap \varphi^{-1}(A) = \emptyset ) = \mathrm{E} \left( \mathrm{E} \left[ \mathbf{1}(X \cap \varphi^{-1}(A) = \emptyset) | M \right] \right) = \\
	& = \mathrm{E} \exp \{ - M(\varphi^{-1}(A)) \} = \mathrm{E} \exp \left\{ - \int_{\varphi^{-1}(A)} \Lambda(u) \, \mathrm{d} u \right\},
\end{align*}
where $\varphi^{-1}(u)=(\delta\|u\|)^{-1}(e^{\|u\|}-1) u$ for $\varphi$ given by \eqref{eq:phi}.
Applying the substitution theorem, we get
\begin{align*}
  \int_{\varphi^{-1}(A)} \Lambda(u) \, \mathrm{d} u = \int_{A} \Lambda(\varphi^{-1}(u)) | D_{\varphi^{-1}}(u) | \, \mathrm{d} u,
\end{align*}
where $D_{\varphi^{-1}}(u)$ is the Jacobian. Hence
\begin{align*}
  \mathrm{P} ( \varphi(X) \cap A = \emptyset ) = \mathrm{E} \exp \left\{ - \int_{A} \Lambda(\varphi^{-1}(u)) | D_{\varphi^{-1}}(u) | \, \mathrm{d} u \right\}.
\end{align*}
Since the distribution of a locally finite simple point process is uniquely determined by the void probabilities, it follows that $\varphi(X)$ is a Cox process with random intensity function
\begin{align*}
  \widetilde{\Lambda}(u) & = \Lambda(\varphi^{-1}(u)) | D_{\varphi^{-1}}(u) | \\ & = \exp \{ \log | D_{\varphi^{-1}}(u) | + \zeta(\varphi^{-1}(u)) + Y_0(\varphi^{-1}(u)) \} = \\
	& = \exp\{ \log | D_{\varphi^{-1}}(u) | + \beta_0 + \beta_1 \| {\varphi^{-1}}(u) \| - \sigma^2/2 + \sigma Q (\varphi (\varphi^{-1}(u))) \} = \\
	& = \exp\{ \log | D_{\varphi^{-1}}(u) | + \beta_0 + \beta_1 \| {\varphi^{-1}}(u) \| - \sigma^2/2 + \sigma Q (u) \}.
\end{align*}
The argument of the exponential above is the sum of the stationary and isotropic Gaussian random field $\sigma Q$ and the deterministic trend given by the other terms, i.e., it is a Gaussian random field with non-constant mean but stationary and isotropic covariance structure. Hence $\varphi(X)$ is a second-order intensity-reweighted stationary LGCP \citep{baddeley:moeller:waagepetersen:00} with a tractable form of the intensity function and the pair-correlation function.

\section{Comparison of the 2- and 3-steps approaches}\label{app:comparison}

Table \ref{t:compar-2-3-steps} provides a comparison between the two step procedure from Section~\ref{subsec:CLE} and our three step procedure from Section~\ref{sec:3-steps} in the context of the simulation study performed for the particle dispersal application in Section~\ref{sec:disp.res}. The comparison is based on 
the root mean squared logarithmic error (RMSLE): For a series of estimates $\{\hat\theta_i:i=1,\ldots,n\}$ of $\theta$, the RMSLE is equal to $\sqrt{(1/n)\sum_{i=1}^n(\log\hat\theta_i-\log\theta)^2}$, and as in Figure~\ref{fig:particles.simul}, in order to obtain robust assessments, each RMSLE was computed by removing the lower and upper 5\% estimates. Here, the logarithmic transformation is used to account for the asymmetric distribution of the estimates of $\sigma$, $1/\alpha$, and $\delta$.

In both procedures, the first step is the same for estimating the first-order parameters $\beta_0$ and $\beta_1$ (hence, the RMSLE ratio is equal to 1). For  estimation of the second-order parameters $\sigma$, $\alpha$, and $\delta$, we notice that the three-steps procedure is more accurate for $\sigma$, overall for $1/\alpha$, and to some extent also for $\delta$.

\begin{table}[ht]
\centering
{\small
\begin{tabular}{cccccc}
  \hline
Simulation   &  \multicolumn{5}{c}{RMSLE ratio}\\
series  &  $2\pi e^{\beta_0}/\beta_1^2$ & $-1/\beta_1$ & $\sigma$  & $1/\alpha$ & $\delta$ \\
\hline
1& 1.00&1.00&    0.18   &   0.73   &   0.71 \\
2& 1.00&1.00&     0.12   &   0.85   &   0.75 \\ 
3& 1.00&1.00&     0.19   &   0.89   &   0.47  \\
4& 1.00&1.00&     0.08   &   0.76   &   1.11  \\
5& 1.00&1.00&     0.23   &   0.86   &   1.18  \\
6& 1.00&1.00&     0.27   &   0.42   &   0.58  \\
7& 1.00&1.00&     0.20   &   1.18   &   1.42  \\
8& 1.00&1.00&     0.16   &   0.78   &   0.78  \\
   \hline
\end{tabular}
}
\caption{\label{t:compar-2-3-steps} Relative estimation accuracy obtained for each parameter using the two or three step procedure. The simulation series are those performed for the particle dispersal study (see Section \ref{sec:disp.res}). The relative estimation accuracy was computed as RMSLE(3-steps procedure)/RMSLE(2-steps procedure), i.e., the ratio of the root mean squared logarithmic errors when using the 3-steps and 2-steps procedure, respectively.}
\end{table}

\subsubsection*{Acknowledgements} 

This work was supported by The Danish Council for Independent Research -- Natural Sciences, grant DFF -- 701400074 `Statistics for point processes in space and beyond', by the `Centre for Stochastic Geometry and Advanced Bioimaging', funded by grant 8721 from the Villum Foundation, by Grant Agency of Czech Republic (Project No.\ 19-04412S). 

\bibliographystyle{royal}
\bibliography{masterbibLGCP}

\end{document}